\documentclass[11pt] {article}      %
\usepackage{amsmath,amssymb, amsthm, eucal,bm}  
\usepackage{graphicx, amscd}
\usepackage{framed}

\let\OLDthebibliography\thebibliography
\renewcommand\thebibliography[1]{
  \OLDthebibliography{#1}
  \setlength{\parskip}{3pt}
  \setlength{\itemsep}{0pt plus 0.3ex}
}

\usepackage{txfonts}                  
\usepackage[T1]{fontenc}     

\usepackage{scrextend}     

\usepackage{enumerate}       
\usepackage{float}                 
\usepackage{tikz}
\usetikzlibrary{matrix, arrows, calc}

\theoremstyle{plain}
\newtheorem{theorem}{Theorem}[section]            
\newtheorem{proposition}[theorem]{Proposition}  

\theoremstyle{definition}
\newtheorem{definition}[theorem]{Definition}

\numberwithin{theorem}{section}
\numberwithin{equation}{section}
\newcommand{\gaction}[2]{\genfrac{}{}{0.5pt}{}{#1}{#2}%
                        \!\lower2pt\hbox{\rotatebox[origin=c]{-90}{{$\looparrowright$}}}}
\newcommand{\dotfraction}[2]{\genfrac{}{}{0.5pt}{}{#1}{#2}%
                        \!\lower.5pt\hbox{{$\circ$}}}

\usepackage{titlesec}                                                           
\titlelabel{\thetitle.\ }                                                         
\titleformat*{\section}{\fontsize{14pt}{14pt} \bf}        
\usepackage{fancyhdr}
\usepackage{sidecap}  
\usepackage[hang,small,sc]{caption}

\def\SL{\hbox{SL}}
\def\vvec#1#2{
{\renewcommand*{\arraystretch}{.7}
\begin{bmatrix}#1\\#2\end{bmatrix}}
}

\def\-{\hbox{\raisebox{.75pt}{-}}}
\def\z{\hbox{\raisebox{.75pt}{-}}}

\def\Euc{{\rm E}}

\makeatletter
\newcommand*\bigcdot{\mathpalette\bigcdot@{.41}}
\newcommand*\bigcdot@[2]{\mathbin{\vcenter{\hbox{\scalebox{#2}{$\m@th#1\bullet$}}}}}
\makeatother
\makeatletter
\newcommand*\bbigcdot{\mathpalette\bigcdot@{.61}}
\newcommand*\bbigcdot@[2]{\mathbin{\vcenter{\hbox{\scalebox{#2}{$\m@th#1\bullet$}}}}}
\makeatother

\def\sz{.02}  

\begin{document}

\title{\bf Apollonian coronas and a new Zeta function}

\author{Jerzy Kocik 
    \\ \small Department of Mathematics, Southern Illinois University, Carbondale, IL62901
   \\ \small jkocik{@}siu.edu
}


\date{\small\today}

\maketitle

\begin{abstract}
\noindent
We  find a formula for the area of disks tangent to a given disk in an Apollonian disk packing 
(corona)  in terms of a certain novel arithmetic Zeta function.
The idea is based on ``tangency spinors'' 
defined for pairs of tangent disks.
\end{abstract}

\section{Introduction}

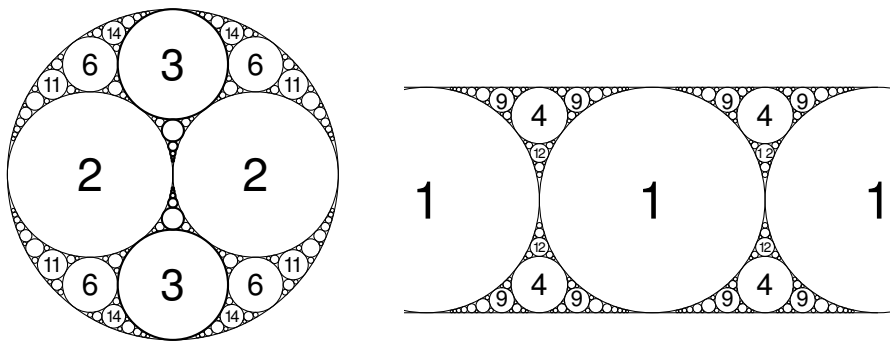
\begin{figure}[H]
\centering
\begin{tikzpicture}[scale=2.2]

\draw (0,0) circle (1);

\foreach \a/\b/\c   in {
1 / 0 / 2 
}
\draw (\a/\c,\b/\c) circle (1/\c)
          (-\a/\c,\b/\c) circle (1/\c);

\foreach \a/\b/\c in {
0 / 2 / 3 ,
0 /4 /15 ,
0 / 6 / 35 , 
0 / 8/ 63,
0 /10 / 99,
0 / 12 / 143
}
\draw[thick] (\a/\c,\b/\c) circle (1/\c)
          (\a/\c,-\b/\c) circle (1/\c) ;

\foreach \a/\b/\c/\d in {
3 / 4 /6, 	8 / 6 / 11,	5 / 12/ 14,	15/ 8 / 18,	8 / 12 / 23,	7 / 24 / 26,
24/	10/	27, 	21/	20/	30, 	16/	30/	35, 	3/	12/	38, 	35/	12/	38, 	24/	20/	39, 	9/	40/	42,
16/	36/	47, 	15/	24/	50, 	 48/	14/	51, 	45/	28/	54, 	24/	30/	59, 	40/	42/	59, 	11/	60/	62,
21/	36/	62, 	48/	28/	63, 	33/	56/	66, 	63/	16/	66, 	8/	24/	71, 	55/	48/	74, 	24/	70/	75,
48/	42/	83, 	80/	18/	83, 	13/	84/	86, 	77/	36/	86, 	24/	76/	87, 	24/	40/	87, 	39/	80/	90
}
\draw (\a/\c,\b/\c) circle (1/\c)       (-\a/\c,\b/\c) circle (1/\c)
          (\a/\c,-\b/\c) circle (1/\c)       (-\a/\c,-\b/\c) circle (1/\c) 
;
\node at (-1/2,0) [scale=1.7, color=black] {\sf 2};
\node at (1/2,0) [scale=1.7, color=black] {\sf 2};
\node at (0,2/3) [scale=1.6, color=black] {\sf 3};
\node at (0,-2/3) [scale=1.6, color=black] {\sf 3};
\node at (1/2,2/3) [scale=1.1, color=black] {\sf 6};
\node at (-1/2,2/3) [scale=1.1, color=black] {\sf 6};
\node at (1/2,-2/3) [scale=1.1, color=black] {\sf 6};
\node at (-1/2,-2/3) [scale=1.1, color=black] {\sf 6};
\node at (8/11,6/11) [scale=.77, color=black] {\sf 1$\!$1};
\node at (-8/11,6/11) [scale=.77, color=black] {\sf 1\!1};
\node at (8/11,-6/11) [scale=.77, color=black] {\sf 1\!1};
\node at (-8/11,-6/11) [scale=.77, color=black] {\sf 1\!1};
\node at (5/14,6/7) [scale=.6, color=black] {\sf 1\!4};
\node at (-5/14,6/7) [scale=.6, color=black] {\sf 1\!4};
\node at (5/14,-6/7) [scale=.6, color=black] {\sf 1\!4};
\node at (-5/14,-6/7) [scale=.6, color=black] {\sf 1\!4};
\end{tikzpicture}
\qquad
%
%
\begin{tikzpicture}[scale=1.5, rotate=90, shift={(0,2cm)}]  
\clip (-1.25,-2.1) rectangle (1.1,2.2);
\draw (1,-2) -- (1,3);
\draw (-1,-2) -- (-1,3);
\draw (0,0) circle (1);
\draw (0,2) circle (1);
\draw (0,-2) circle (1);

\foreach \a/\b/\c in {
3/4/4,  5/12/12,  7/24/24, 9/40/40  
}
\draw (\a/\c,\b/\c) circle (1/\c)    (\a/\c,-\b/\c) circle (1/\c)
          (-\a/\c,\b/\c) circle (1/\c)    (-\a/\c,-\b/\c) circle (1/\c)   ;

\foreach \a/\b/\c in {
8/    6/   9, 	
15/   8/   16 , 	
24/  20/  25, 	
24/  10/  25, 	
21/  20/  28,
16/	30/	33,    
35/	12/	36,
48/	42/	49,
48/	28/	49,
48/	14/	49,
45/	28/	52,
40/	42/	57,
33/	56/	64,
63/	48/	64,
63/	16/	64,
55/	48/	72,
24/	70/	73,
69/	60/	76,
80/	72/	81,
64/	60/	81,
80/	36/	81,
80/	18/	81
}
\draw (\a/\c, \b/\c) circle (1/\c)          (-\a/\c, \b/\c) circle (1/\c)
          (\a/\c,-\b/\c) circle (1/\c)         (-\a/\c,-\b/\c) circle (1/\c)
          (\a/\c,2-\b/\c) circle (1/\c)       (-\a/\c,2-\b/\c) circle (1/\c)
          (\a/\c,2+\b/\c) circle (1/\c)       (-\a/\c,2+\b/\c) circle (1/\c)
          (\a/\c,-2+\b/\c) circle (1/\c)       (-\a/\c,-2+\b/\c) circle (1/\c)
;
\node at (0,0) [scale=1.9, color=black] {\sf 1};
\node at (0,-2) [scale=1.9, color=black] {\sf 1};
\node at (0,2) [scale=1.9, color=black] {\sf 1};
\node at (-3/4,1) [scale=1.1, color=black] {\sf 4};
\node at (-3/4,-1) [scale=1.1, color=black] {\sf 4};
\node at (3/4,1) [scale=1.1, color=black] {\sf 4};
\node at (3/4,-1) [scale=1.1, color=black] {\sf 4};
\node at (7/8,8/6) [scale=.8, color=black] {\sf 9};
\node at (7/8,-8/6) [scale=.8, color=black] {\sf 9};
\node at (7/8,4/6) [scale=.8, color=black] {\sf 9};
\node at (7/8,-4/6) [scale=.8, color=black] {\sf 9};
\node at (-7/8,8/6) [scale=.8, color=black] {\sf 9};
\node at (-7/8,-8/6) [scale=.8, color=black] {\sf 9};
\node at (-7/8,4/6) [scale=.8, color=black] {\sf 9};
\node at (-7/8,-4/6) [scale=.8, color=black] {\sf 9};
\node at (-0.42,1) [scale=.5, color=black] {\sf 1\!2};
\node at (-0.42,-1) [scale=.5, color=black] {\sf 1\!2};
\node at (.42,1) [scale=.5, color=black] {\sf 1\!2};
\node at (.42,-1) [scale=.5, color=black] {\sf 1\! 2};
\end{tikzpicture}
\caption{Apollonian Window (left) and Apollonian Belt (right)}
\label{fig:Apollo}
\end{figure}

{\bf Apollonian disk packing} is an arrangement of infinite number of disks determined by 
three mutually tangent disks through  recurrent inscribing new disks in triangle-like spaces between disks  
as they appear in this process. 
An Apollonian disk packing is called {\bf integral} if the curvatures of all disks are integers.
If any four mutually tangent disks (called disks in {\bf Descartes configuration}) are integral,
so is the whole packing. 
See Figure \ref{fig:Apollo} for two basic examples, the Apollonian Window \cite{jk-t} and the Apollonian Belt.

\noindent
Note that in both arrangements, the disks -- if tangent --  are tangent externally.
The greatest circle in the Apollonian Window is understood as the boundary of the disk that lies outside of it; 
consequently its radius and curvature is equal to $-1$ rather than $1$.
Similarly, the lines in the Apollonian Belt are boundaries if half-planes that lie outside of the central belt and 
are viewed as disks of curvature 0.
The integral packings are interesting due to their evident connection with number theory  
\cite{Hir,GLM1,LMW,Man,Mel}.
\\

\begin{figure}
\begin{center}
\begin{tikzpicture}[scale=3, rotate=90]  
\clip (-1.2,-.2) rectangle (.2,2.2);


\draw (1,-1) -- (1,3);
\draw (-1,-1) -- (-1,3);

\draw (0,0) circle (1);
\draw (0,2) circle (1);


\foreach \a/\b/\c in {
3/  4/ 4
}
\draw 
          (-\a/\c,\b/\c) circle (1/\c)  
;

\foreach \a/\b/\c in {
8/    6/   9,
15/   8/   16 ,
24/  20/  25,
24/  10/  25,
35/	12/	36,
48/	42/	49,
48/	28/	49,
48/	14/	49,
63/	48/	64,
63/	16/	64,
80/	72/	81,
80/	36/	81,
80/	18/	81,
99/	60/	100,
99/	20/	100,
120/	110/	121,
120/	88/	121,
120/	66/	121,
120/	44/	121,
120/	22/	121,
143/	120/	144,
143/	24/	144,
168/	156/	169,
168/	130/	169,
168/	104/	169,
168/	78/	169,
168/	52/	169,
168/	26/	169,
189/	148/	196,
195/	140/	196,
195/	84/	196,
195/	28/	196
}
\draw 
          (-\a/\c,\b/\c) circle  (1/\c) 
          (-\a/\c,2-\b/\c) circle (1/\c)
;
\draw
node at (-1.1,2) {$\frac{0}{1}$}
node at (-1.1,1.33) {$\frac{2}{3}$}
node at (-1.1,1) {$\frac{1}{1}$}
node at (-1.1,0.66) {$\frac{4}{3}$}
node at (-1.1,0) {$\frac{2}{1}$}
;
\end{tikzpicture}
\end{center}
\caption{Ford circles}
\label{fig:ford}
\end{figure}
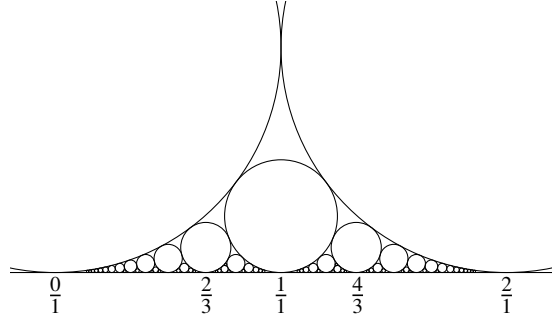

Let us define a {\bf corona} of a particular disk
in an Apollonian disk packing $\mathcal A$ as the set of all disk that are tangent to $\mathcal A$.
Our point of interest is the area of the  coronas or its fragments.
The famous Ford circles \cite{Ford} can be interpreted as a corona in the Apollonian Belt 
composed of circles tangent to 0-curvature disk at the bottom (see Figure \ref{fig:ford}).
The area of the disks between two large circles  was found in \cite{WB}.  
Quite remarkably, it is related to Riemann zeta functions:
\begin{equation}
\label{eq:Ford}
A_{\rm Ford} 
 \ = \  \pi \left( \, 1+ \frac{\zeta(3)}{\zeta(4)} \right)
\end{equation}

~

In this paper we present a formula for coronas and their fragments in an arbitrary Apollonian disk packing.
To give a taste of its arithmetic elegance, here is a specific result: 
 the corona made of the disks tangent to the 
greatest circle (of radius 1, see Figure \ref{fig:Apollo}, left) of the Apollonian Window
consists of four congruent quarters.
The area of each is:
\begin{equation}
\label{eq:nice}
A_{\rm AW} \ = \ \pi\!\!\!\!\! \sum_{\substack{0\,\leq \,k\leq \,n  \\ \gcd(n,k)=1}}\ \frac{1}{\big( n^2+k^2+1\big)^2}
\end{equation}
where the sum is restricted to coprimes $(n,k)$.
The greatest common divisor will be denoted $(m,n)$. 
One easily recognizes the largest disks among the first summands: 
$$
A_{\rm AW} \ = \  \frac{1}{\big( 1^2 \!+\! 0^2 \!+\! 1\big)^2}  \!+\!     
                      \frac{1}{\big(1^2 \!+\! 1^2 \!+\! 1\big)^2}  \!+\! 
                      \frac{1}{\big(2^2 \!+\! 1^2 \!+\! 1\big)^2}  \!+ ...
\ = \ 
\  \frac{1}{2^2} + \frac{1}{ 3^2} + \frac{1}{6^2} + ...
$$
Section \ref{sec:corona} provides the general formula for coronas in arbitrary Apollonian disk packings. 
The formula defines a new type of arithmetic function quite similar to Epstein Zeta function,
except for a free term (1 in the above formula) 
and, more importantly,  the sum being taken over coprimes instead of all pairs of integers.

The idea is based on the existence of the ``spinor structure'' in the Apollonian disk packing.

\section{Spinors in Apollonian disk arrangements}

The section summaries the essential facts extracted from \cite{jk-t}.
\\

\begin{figure}[h]
\centering
\includegraphics[scale=.8]{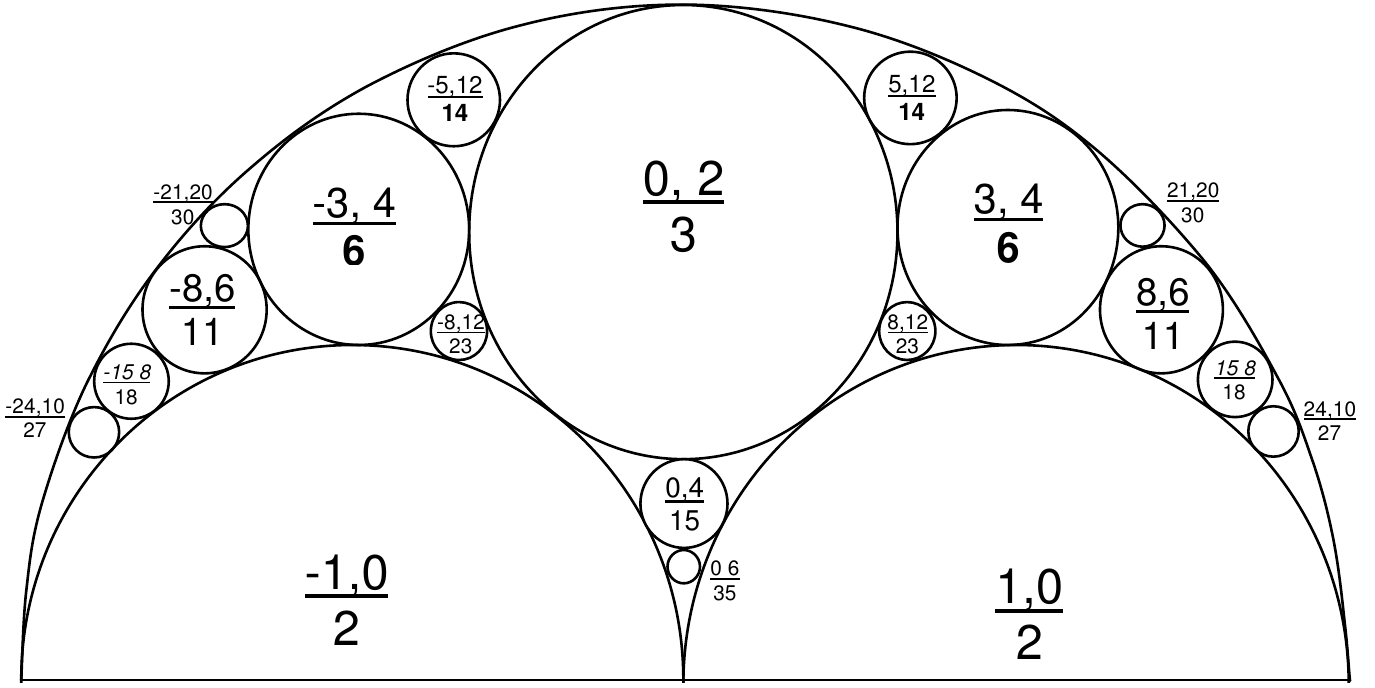}
\caption{\small Symbols in the Apollonian Window (upper half).}
\label{fig:S32}
\end{figure}

Every disk in the Cartesian plane may be given a {\bf symbol},
a fraction-like label that encodes the size and position of the disk:  
the curvature (reciprocal of radii) is indicated in the denominator 
while the positions of the centers may be read off by interpreting the symbol as a pair of fractions
\cite{jk-c}.
For example:
$$
\hbox{\sf symbol:\  } 
       \frac{3,\;4}{6}  
                         \qquad \Longrightarrow \qquad 
       \begin{cases}     \hbox{\sf radius:\  } &  r = \frac{1}{6} \\
                                  \hbox{\sf center:\ }  & (x,\,y)  =   \left(\frac{3}{6},\, \frac{4}{6}\right)
                                                                    =   \left(\frac{1}{2},\, \frac{2}{3}\right)
       \end{cases}
$$
The numerator will be called the {\bf reduced coordinates} of the a disk's center 
and denoted by dotted letters $(\dot x, \, \dot y)= (x/r,\, y/r)$.
Unbounded disk extending outside a circle are given negative radius and curvature.
Quite remarkably, all symbols in the Apollonian Window have integer entries
(see Figure~\ref{fig:S32}).

The symbols in an Apollonian disk packing can be generated from the first three via the Descartes theorem. 
It states that the curvatures of four mutually tangent disks, the so-called {\bf Descartes configuration}, 
satisfy the Descartes formula (1643)  \cite{Coxeter,Sod}:
\begin{equation}
\label{eq:Descartes}
                  (A+B+C+D)^2 = 2\; (A^2 + B^2 + C^2 + D^2 )
\end{equation}
A more convenient version of this quadratic formula is
\begin{equation}
\label{eq:Boyd}
                 D+D' = 2(A+B+C)
\end{equation}
where  $D$ and $D'$ are the two solutions to \eqref{eq:Descartes}, given $A$, $B$, $C$.
The last formula holds for the corresponding reduced positions of the centers of the disks  
(consult \cite{jk-d}).
Since an Apollonian gasket is a completion of a Descartes configuration, 
all symbols may be derived from the first four using (\ref{eq:Boyd}) . 
Thus integrality of the symbols in the Apollonian Window follows from the integrality of the first four 
disks.
%
\begin{figure}[h]
\centering
\includegraphics[scale=.78]{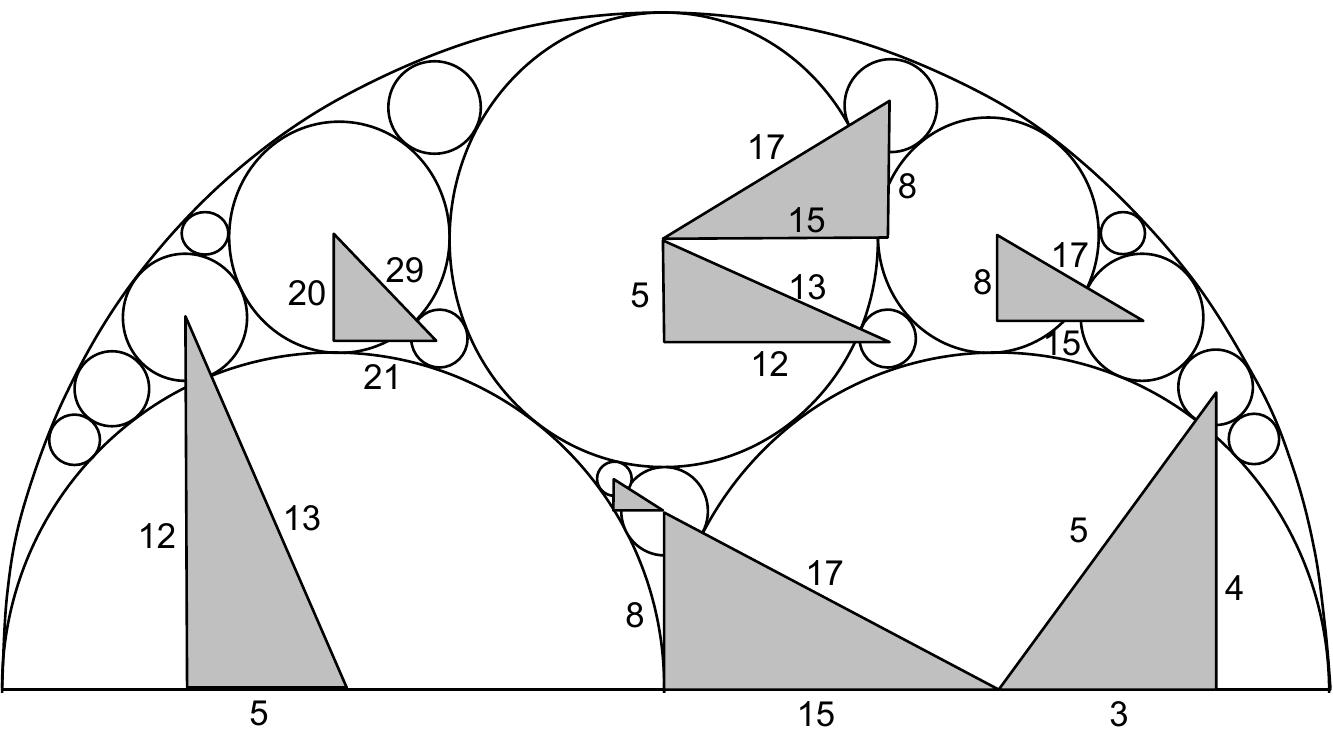}
\caption{\small Pythagorean triples in the Apollonian Window}
\label{fig:AWtriangles}
\end{figure}


Every pair of tangent disks in $\mathcal A$ defines a right triangle
with sides proportional to a {\bf Pythagorean triple} as follows:
\begin{equation}
\label{eq:product}
     \frac{\dot x_1, \;  \dot y_1}{ \beta_1}
\ \Join  \
     \frac{\dot x_2,\; \dot y_2}{ \beta_2}
\qquad \mapsto \qquad
       \left[\begin{matrix}
                     a \\
                     b \\
                    c
     \end{matrix}\right]
\equiv 
  \left[\begin{matrix}
                     \beta_1 \dot x_2 - \beta_2 \dot x_1 \\
                     \beta_1 \dot y_2 - \beta_2 \dot y_1 \\
                    \beta_1+\beta_2
     \end{matrix}\right]
\end{equation}
where $a^2+b^2=c^2$ 
(see Figure \ref{fig:AWtriangles} for examples).
Note that in the case of the Apollonian Window, all of these Pythagorean triples are integral.
The actual size of the triangle in the plane is scaled down by the factor of $\beta_1\beta_2$
(gray triangles in Figure \ref{fig:AWtriangles}).

The next step is to recall that Pythagorean triples admit Euclidean parameters 
that determine them via the following prescription:
\begin{equation}
\label{eq:euclid}
u=[m,n] \quad \to \quad (a,b,c) = (m^2-n^2,\ 2mn, \ m^2+n^2)
\end{equation}
(see, e.g., \cite{Sie, T-T}).
This map has a simple form in terms of complex numbers:
\begin{equation}
\label{eq:square}
            u=  
            m+ni \quad \to \quad u^2  = a+bi   = (m^2-n^2) + 2mn\, i 
\end{equation}
with $c=|u^2|=m^2+n^2$.  
As explained in  \cite{jk-c}, Euclidean parameters can be viewed as a {\it spinor},
a vector $\mathbf u\in \mathbb R^2$.
Hence our concept of {\bf tagency spinor} defined for an ordered pair of  two tangent disks 
(not necessarily integral)
of curvatures $A$ and $B$ as
as :
\begin{equation}
\label{eq:defspinor}
u = \pm\sqrt{\frac{z}{AB}}
\end{equation}
where $z\in \mathbb C$ is the complex number reprezenting vector joining the centers of the disks.
We view spinor as a vector in a two dimensional Euclidean space and use 
the complex structure as a convenient description,  $u\in\mathbb C\cong\mathbb R^2$
via identification 
$$
\vvec{m}{n} \ \equiv \  m+ni
$$ 
\\

\noindent
{\bf Remark:} 
A Pythagorean triangle with sides satisfying $a^2+b^2-c^2=0$ may be viewed as 
an isotropic vector of Minkowski space $\mathbb R^{2,1}$ and,
as such, it may be represented as the tensor square of a spinor 
from the associated two-dimensional spinor space $\mathbb R^2$.

~~

%
\begin{SCfigure}[2.2][h]
\centering
\includegraphics[width=0.45\textwidth]{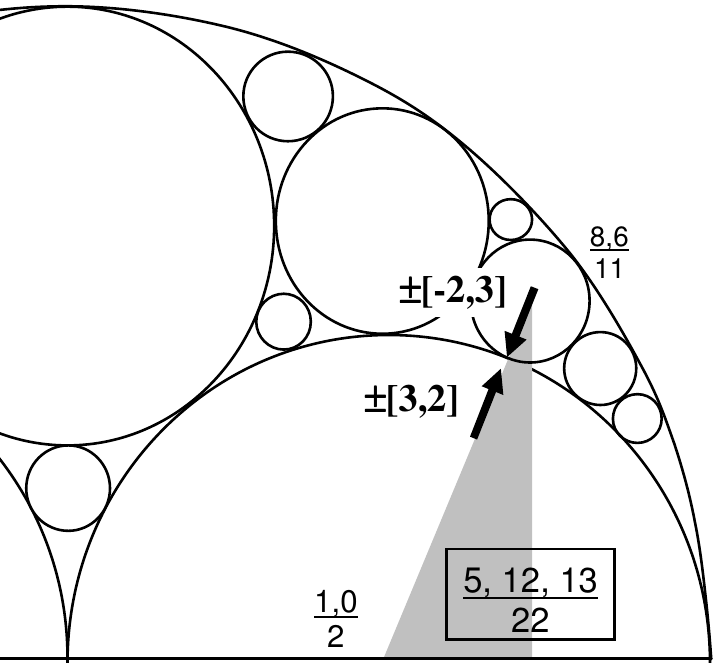}
\caption{\small Spinors in the Apollonian Window}
\label{fig:21}
\end{SCfigure}

~~

%
%

In graphical representation we shall mark a spinor by an arrow that indicates 
the order of disks, and will label it by its value.
Note that the spinor is defined up to a sign, since $(-u)^2=u^2$.
Also, the spinor depends on the order of the circles: 
 if $u$ is a spinor for $(AB)$,  then the spinor for  $(BA)$  is $u' = iu$  (up to a sign).
\\

%
%

\begin{figure}[H]
\centering
\includegraphics[scale=.7]{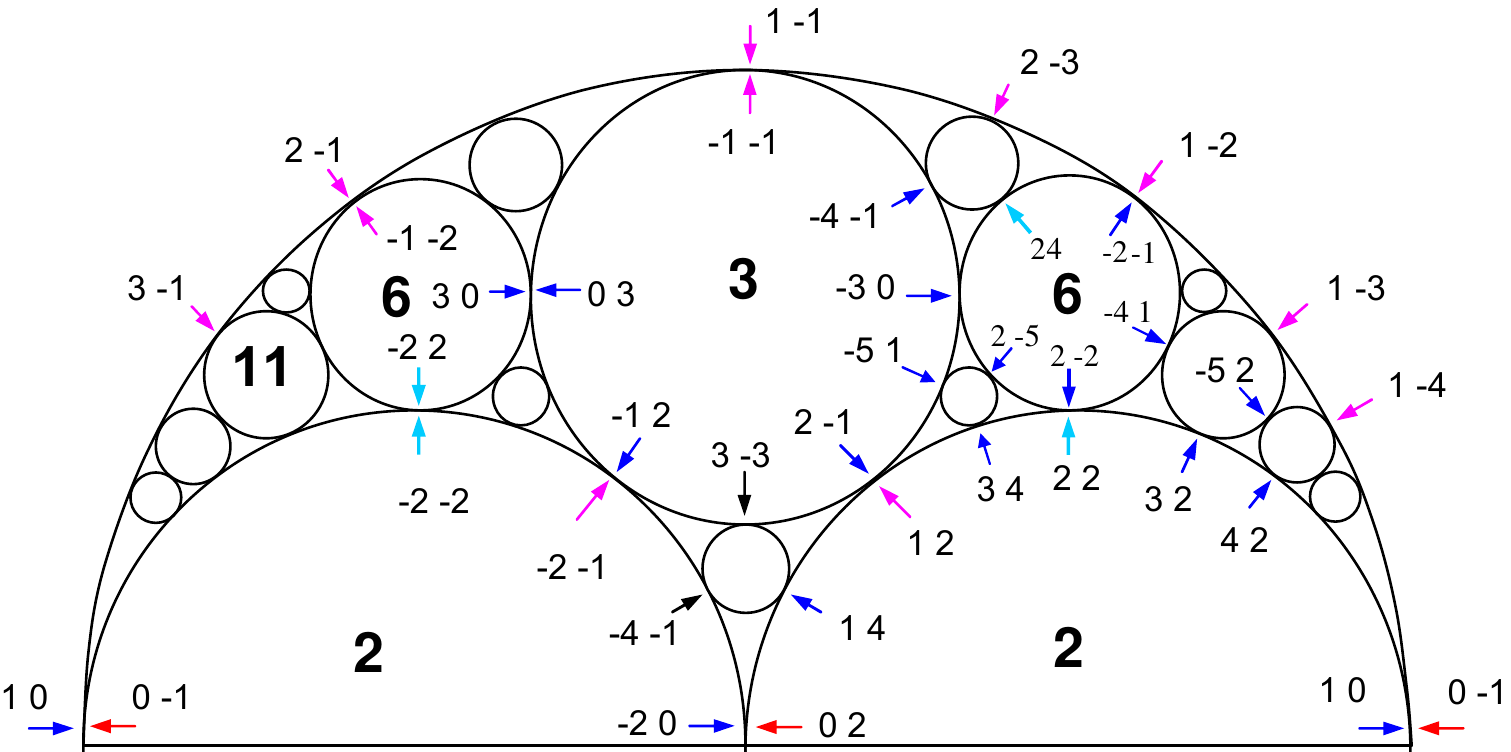}
\caption{\small Flow of spinors in the Apollonian Window}
\label{fig:AWspinors}
\end{figure}

Figure \ref{fig:AWspinors} shows spinors in the Apollonian Window $\mathcal A$   
(for visual clarity the brackets are omitted). 
There are four vectors at every point of tangency:  two signs times two orderings of circles. 
All spinors in the Apollonian Window $\mathcal A$ are  integral but
we will not restrict our considerations to the integral examples.  
%
%


The key property of spinors is that in a quartet of mutually tangent disks
(disks in a Descartes configuration) 
they admit a choice of signs $(\pm)$ 
such that these two properties hold
$$
\begin{array}{cc}
    ``\hbox{curl}\; \mathbf u = 0\hbox{''} :   &\qquad   u_{12} + u_{23} + u_{31} = 0\\
    ``\hbox{div}\; \mathbf u = 0\hbox{''} :   &\qquad   u_{14} + u_{24} + u_{34} = 0\\
\end{array}
$$
where $u_{ij}$ represents a spinor for  $i$-th and $j$-th tangent disks.
These results may be viewed as am ``underground'' version of Descartes' theorem on circles.


These laws have local character.  Extending the appropriate choice of signs to a greater system, 
like the whole Apollonian gasket,  will encounter topological obstructions.

~\\


~\\
In the remainder of this section we review the main properties of tangency spinors; for proofs see \cite{jk-t}.
The capital letters will denote both circles and their curvatures.
\\

%
%
%

\begin{figure}[h!]
\centering
\includegraphics[scale=.8]{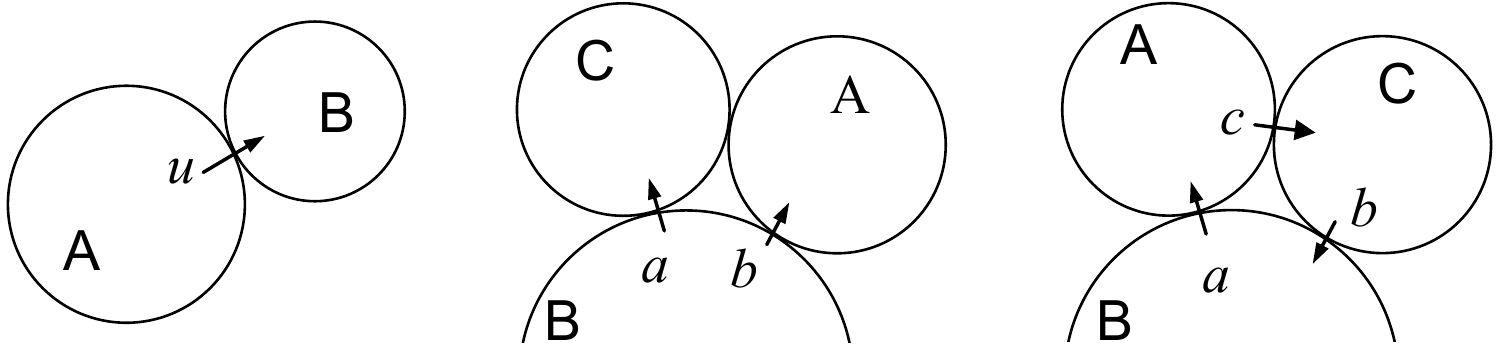}
\caption{\small Left: Two curvatures and a spinor; Center: curvature from two spinors; Right: curl u = 0  }
\label{fig:uuu}
\end{figure}

\begin{proposition} 
\label{thm:curvs}
 If $u$ is the tangency spinor for two tangent disks of curvatures $A$ and $B$, respectively,   
(Figure \ref{fig:uuu}, left)  then 
\begin{equation}
\label{eq:thm3}
                     |u|^2 = A + B
\end{equation}
\end{proposition}


\begin{theorem}[\bf curvatures from spinors]    
\label{thm:curv}
In the system of three mutually tangent circles, the symplectic product of two spinors directed outward from 
(respectively inward into) one of the circles equals (up to sign) its curvature, 
e.g., following notation of Figure \ref{fig:uuu}, center:
\begin{equation}
\label{eq:thm4}
          B  \  =  \ \pm\;  a\times b
\end{equation}
where  $ a\times b\   :=   \  \det [a|b]  \ = \ a_1b_2-a_2b_1$.
\end{theorem}

%
%

\begin{theorem}[\bf spinor curl] 
\label{thm:curl}
The signs of the three tangency spinors between be three mutually tangent circles
(Figure \ref{fig:uuu}, right)  
may be chosen so that 
\begin{equation}
\label{eq:curl}
a  +  b  +  c  =  0  \qquad\quad         [ \hbox{\rm ``curl}\ \mathbf u = 0 \hbox{''}]   
\end{equation}
\end{theorem}

~\\

\begin{figure}[h]
\centering
\includegraphics[scale=.8]{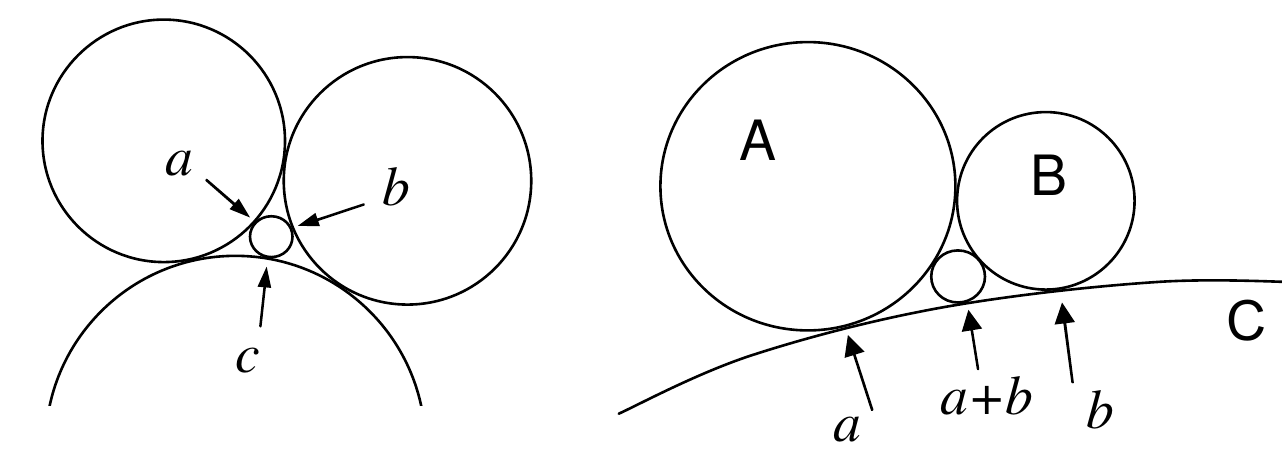}
\caption{\small (a) vanishing divergence, (b) spinor addition}
\label{fig:thmdiv}
\end{figure}

\begin{theorem}  
\label{thm:div}
Let $A$, $B$, $C$, and $D$ be four circles in a Descartes configuration.  
\\[5pt]
{\bf [Vanishing divergence]:}  
If $a$, $b$ and $c$ are tangency spinors for pairs $AD$, $BD$ and $CD$  
(Figure \ref{fig:thmdiv} left),
then their signs may be chosen so that 
\begin{equation}
\label{eq:thm5a}
 					a + b + c = 0 	\qquad	[\hbox{``\rm div}\, \mathbf u = 0\hbox{''}]                                   
\end{equation}
The same property holds for the outward oriented spinors.
\\[7pt]
{\bf [Additivity]:}  If $a$ and $b$ are tangency spinors for pairs $CA$ and $CB$ 
( Figure \ref{fig:thmdiv} right), 
then there is a choice of signs such that the sum 
\begin{equation}
\label{eq:thm5b}
                                      c = a + b
\end{equation}
is a spinor of tangency for $CD$.  
\end{theorem}

%
%
%

Theorem \ref{thm:div} can be iterated to produce spinors for all circles 
inscribed between two initial circles. 
\\

The last m concerns the signs of the spinors.
Let $A$, $B$ and $C$ be three mutually tangent disks,
and let $D$ and $D'$ be two disks that complete this triple to respective Descartes configurations.
Consider spinors $a$, $b$ and $d$  from circle $C$ 
to circles $A$, $B$, and $D$,
like in the Figure \ref{fig:thmdiv}, with some fixed singes. 
$D$ is the unlabeled small disk in the center.
Spinors $a$ to $b$ are {\bf harmonized} over the arc of circle $C$ through $D$ if
$a+b=\pm d$ for some choice of the sign of $d$.
If $a$ and $b$ are harmonized, so are $(-a)$ and $(-b)$. 
If  (signs of ) spinors are not harmonized, then spinors
$a$ and $(-b$ are harmonized over the complementary part of the circle $C$,.i.e., the arc tangent to $D'$.
\\

Here is an example referring to Figure \ref{fig:AWspinors}:
%
%
%
Spinors 
$\vvec12$ and $\vvec22$ 
are harmonized over the upper arc of disk ``2''
since $\vvec12+\vvec22=\vec34$. 
But  
$\vvec12$ and $\vvec{-2}{-2}$ are not, they are harmonized over the lower arc since
$\vvec12+\vvec{-2}{-2}=-\vvec10$,
a spinor towards the external disk ``$-1$, which is the other solution to Descartes problem for disks $(2,3,6)$.
\\

Now we move to our main theorem on the areas of Apollonian coronas.

\newpage

\section{Corona's area}
\label{sec:corona}

Here is the main result:

\begin{theorem}
\label{thm:f}
Let $\mathbf u$ and $\mathbf v$ be two spinors oriented from a disk of curvature $B$ 
in an Apollonian disk packing 
towards  two mutually tangent disk in its corona. Define  matrix $M$,
and a ``dummy'' integer vector $\mathbf f$ as:
$$
\mathbf u=\begin{bmatrix} a \\ b  \end{bmatrix}, \ \  
\mathbf v=\begin{bmatrix} c \\ d  \end{bmatrix}, \quad 
\qquad M=\begin{bmatrix} a & c  \\ b & d \end{bmatrix}, 
\quad
\mathbf f = \begin{bmatrix}m\\ n \end{bmatrix}
$$
Then the area of the entire corona is
\begin{equation}
\label{eq:main}
\begin{aligned} 
         A(M)  &= \ \frac{\pi}{2} \ \sum_{\mathbf f \,\in\, \mathbb Z^2_{\rm o}} \ 
                        \frac{1}{\big( \|\; M\, \mathbf f \;\|^2  - B\big)^{-2}} \\
                  & = \ \frac{\pi}{2} \ \sum_{\substack{m,n\,\in \mathbb Z \\ (m,n)=1}} \ 
                                             \frac{1}{ \big( (am+cn)^2 + (bm+dn)^2 - B\big)^2}
\end{aligned}
\end{equation}
\end{theorem}

\noindent
{\bf Towards the proof.}
Suppose we want to find the area of a fragment of corona
consisting of disks between the circles $A$ and $B$
around the base circle of curvature $C$ in Figure \ref{fig:thmdiv}.
We can recover all spinors between $\mathbf a$ and $\mathbf b$ using recursively Theorem \ref{thm:div}B.  
in a form of a series of sequences $R_i$, each obtained from the previous by inserting 
new terms that are sums of the neighboring terms of the previous. 
A few initial rows are shown below:
$$
\begin{array}{ccccccccccccccccc}
R_0: \
& \mathbf a  &         &        &            &        &        &          &            &\mathbf b \\[11pt] 
R_1:\
&\mathbf a  &         &        &            &\mathbf a+\mathbf b   &        &         &             &\mathbf b \\[11pt] 
R_2:\
& \mathbf a  &         &2\mathbf a+\mathbf b &           &\mathbf a+\mathbf b  &         & \mathbf a+2\mathbf b  &            &\mathbf b \\[11pt] 
R_3:\
& \mathbf a  &3\mathbf a+\mathbf b  &2\mathbf a+\mathbf b &3\mathbf a+ 2\mathbf b &\mathbf a+\mathbf b  &2\mathbf a+3\mathbf b & \mathbf a+2\mathbf b  &2\mathbf a+3\mathbf b &\mathbf b \\[11pt] 
\end{array}
$$
(If $R_{ni}$ denotes the i-th term in the n-th row, 
the iteration is defined by $R_{n+1, 2i}=R_{n, i}$   and  $R_{n+1, 2i+1}=R_{n, i}+R_{n,i+1}$,
where $i=0,...,2^n$).
Such a system of sequences will be called a (generalized)  Stern-Brocot array.
After removal of the repetitious terms along the columns, it becomes a Stern-Brocot tree. 

Each of the spinors $\mathbf u$ from this collection determines the curvature of the associated circle
tangent to $C$:
\begin{equation}
\label{eq:use}
curv = \|\mathbf u\|^2 - C
\end{equation}
by Proposition \ref{thm:curvs}, Eq. \eqref{eq:thm3}.
The corresponding area is
$$
area = \frac{\pi}{(curv)^2}
$$
and the total area is the sum over all entries produced by the above scheme
(metaphorically, all entries in the ``last'', limit, row $R_\infty$). 

\def\sz{.25}

\begin{figure}
\begin{tikzpicture}[scale=6.7]
\clip (-.1,-.2) rectangle (1.1,1.1);
\draw (0,0) circle (1);

\draw [fill=yellow, opacity=.35]   (1/2 ,0/2) circle (1/2);
\draw [fill=blue, opacity=.25, thick]   (0/3 ,2/3) circle (1/3);    


\foreach \a/\b/\c/\d in {
3 / 4 /6,          
8 / 6 / 11,         
15/ 8 / 18,         
8 / 12 / 23,
24/	10/	27,         
3/	12/	38,
35/	12/	38,         
24/	20/	39,
15/	24/	50,
48/	14/	51,         
24/	30/	59,
48/	28/	63,
63/	16/	66,         
8/	24/	71,
48/	42/	83,
80/	18/	83,
24/	40/	87,
80/	36/	95,
63/	48/	98,
3/	20/	102,
99/	20/	102,
8/	30/	107,
48/	56/	111,
15/	40/	114,
120/	22/	123,
35/	60/	134,
99/	60/	134,
120/	44/	135,
80/	72/	143,
143/	24/	146,
48/	70/	147,
120/	66/	155,
63/	80/	162,
24/	60/	167,
168/	26/	171,
80/	90/	179,
120/	88/	183,
168/	52/	183,
48/	84/	191,
3/	28/	198,
195/	28/	198,
8/	42/	203
0 / 2 / 3,  
0 /4 /15 , 
0 / 6 / 35, 
0 / 8/ 63,
0 /10 / 99,
0 / 12 / 143,
0 / 14 / 195,
0 / 16 / 255,
0 / 18 / 323
}
\draw  [fill=blue, opacity=\sz, thick] (\a/\c,\b/\c) circle (1/\c)      
          (\a/\c,-\b/\c) circle (1/\c)      
;

\foreach \a/\b/\c/\s   in {
1 / 0 / 2 / 3
}
\draw [thick]  (\a/\c,\b/\c) circle (1/\c)
          (-\a/\c,\b/\c) circle (1/\c)
node  [scale=\s]  at (\a/\c,\b/\c)   {$\c$}
node  [scale=\s]  at (-\a/\c,\b/\c) {$\c$};

\foreach \a/\b/\c/\s in {
0 / 2 / 3   /3,
0 /4 /15   /1,
0 / 6 / 35 / .6 
}
\draw (\a/\c,\b/\c) circle (1/\c)
          (\a/\c,-\b/\c) circle (1/\c)
        node [scale= \s]  at (\a/\c,\b/\c)  {$\c$}
        node [scale= \s]  at (\a/\c,-\b/\c) {$\c$};

\foreach \a/\b/\c in {
0 / 8/ 63,
0 /10 / 99
}
\draw (\a/\c,\b/\c) circle (1/\c)  
           (\a/\c,-\b/\c) circle (1/\c) 
node [scale= 20/\c]  at (\a/\c,\b/\c)  {$\c$}
node [scale= 20/\c]  at (\a/\c,-\b/\c) {$\c$};

\foreach \a/\b/\c in {
0 / 12 / 143,
0 / 14 / 195,
0 / 16 / 255,
0 / 18 / 323
}
\draw (\a/\c,\b/\c) circle (1/\c)  
           (\a/\c,-\b/\c) circle (1/\c);

\foreach \a/\b/\c/\s in {
3 / 4 /6/3,
8 / 6 / 11/1.5,
5 / 12/ 14/1.2
}
\draw (\a/\c,\b/\c) circle (1/\c)       (-\a/\c,\b/\c) circle (1/\c)
          (\a/\c,-\b/\c) circle (1/\c)       (-\a/\c,-\b/\c) circle (1/\c)
node [scale=\s]  at (\a/\c,\b/\c)    {$\c$}
node [scale=\s]  at (-\a/\c,\b/\c)   {$\c$}
node [scale=\s]  at (\a/\c,-\b/\c)   {$\c$}
node [scale=\s]  at (-\a/\c,-\b/\c)  {$\c$};

\foreach \a /  \b / \c   in {
15/ 8 / 18,
8 / 12 / 23,
7 / 24 / 26,
24/	10/	27,
21/	20/	30,
16/	30/	35,
3/	12/	38,
35/	12/	38,
24/	20/	39,
9/	40/	42,
16/	36/	47,
15/	24/	50,
48/	14/	51,
45/	28/	54,
24/	30/	59,
40/	42/	59,
11/	60/	62,
21/	36/	62,
48/	28/	63
}
\draw (\a/\c,\b/\c) circle (1/\c)          (-\a/\c,\b/\c) circle (1/\c)
          (\a/\c,-\b/\c) circle (1/\c)         (-\a/\c,-\b/\c) circle (1/\c)
node [scale=20/\c]  at (\a/\c,\b/\c)    {$\c$}
node [scale=20/\c]  at (-\a/\c,\b/\c)   {$\c$}
node [scale=20/\c]  at (\a/\c,-\b/\c)   {$\c$}
node [scale=20/\c]  at (-\a/\c,-\b/\c)  {$\c$};


\foreach \a /  \b / \c   in {
33/	56/	66,
63/	16/	66,
8/	24/	71,
55/	48/	74,
24/	70/	75,
48/	42/	83,
80/	18/	83,
13/	84/	86,
77/	36/	86,
24/	76/	87,
24/	40/	87,
39/	80/	90,
64/	60/	95,
80/	36/	95,
33/	72/	98,
63/	48/	98,
65/	72/	98,
3/	20/	102,
99/	20/	102,
8/	30/	107,
56/	90/	107,
5/	36/	110,
69/	60/	110,
91/	60/	110,
48/	92/	111,
48/	56/	111,
15/	112/	114,
15/	40/	114,
40/	72/	119,
39/	96/	122,
120/	22/	123,
117/	44/	126,
32/	126/	131,
48/	102/	131,
112/	66/	131,
35/	60/	134,
99/	60/	134,
120/	44/	135,
105/	88/	138,
32/	132/	143,
80/	72/	143,
17/	144/	146,
143/	24/	146,
48/	70/	147,
96/	110/	147,
45/	76/	150,
51/	140/	150,
120/	66/	155,
85/	132/	158,
63/	80/	162,
24/	60/	167,
56/	120/	167,
136/	84/	167,
55/	96/	170,
119/	120/	170,
72/	154/	171,
168/	26/	171,
69/	92/	174,
165/	52/	174,
64/	102/	179,
80/	90/	179,
160/	78/	179,
19/	180/	182,
51/	156/	182,
141/	84/	182,
120/	88/	183,
168/	52/	183,
57/	176/	186,
153/	104/	186,
48/	84/	191,
128/	132/	191,
65/	120/	194,
95/	168/	194,
129/	120/	194,
144/	130/	195,
3/	28/	198,
195/	28/	198,
8/	42/	203
}
\draw (\a/\c,\b/\c) circle (1/\c)          (-\a/\c,\b/\c) circle (1/\c)
          (\a/\c,-\b/\c) circle (1/\c)         (-\a/\c,-\b/\c) circle (1/\c);

\draw [->, line width=3pt, purple]  (.93,0) -- (1.07,0);
\node at (.85, 0) {$\begin{bmatrix} 1\\ 0 \end{bmatrix}$};

\draw [->, line width=3pt, purple] (0.07,0) -- (-0.07,0);
\node at (.17, 0.0) {$\begin{bmatrix} 0\\ 2 \end{bmatrix}$};

\draw [->, line width=3pt, purple] (0.5,.43) -- (0.5,0.57);
\node at (.5, 0.32) {$\begin{bmatrix} 2\\ 2 \end{bmatrix}$};

\draw [->, line width=3pt, purple] (0.24,.36) -- (0.15,0.47);
\node at (.27, 0.27) {$\begin{bmatrix} 1\\ 2 \end{bmatrix}$};

\end{tikzpicture}
\caption{Apollonian Window -- right op corner fragment  the major corona}
\label{fig:quarter2}
\end{figure}
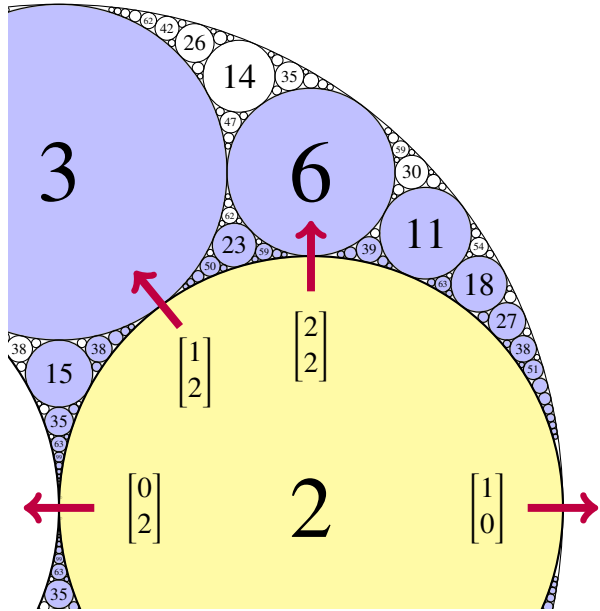

Let us consider a particular example of the corona around the circle of curvature 2 in the Apollonian Window,
shown in Figure \ref{fig:quarter2}.
If we start with vectors $\mathbf a = [1\, 2]^T$ and $\mathbf b = [2\, 2]^T$,  
the Stern-Brocot array looks like the one below, left:
\begin{figure}[H]
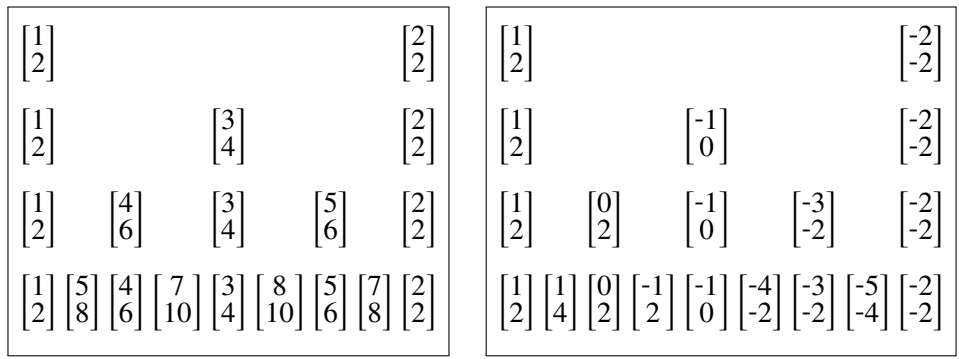

$$
\arraycolsep=1.7pt
\boxed{\begin{array}{ccccccccccccccccc}
~  \\[-10pt]
\vvec 12  &               &             &                  &              &                  &               &             &\vvec 22\\[14pt] 
\vvec 12  &               &             &                  &\vvec 34  &                  &               &             &\vvec 22\\[14pt] 
\vvec 12  &               &\vvec 46 &                  &\vvec 34  &                  & \vvec 56  &             &\vvec 22\\[14pt] 
\vvec 12  &\vvec58    &\vvec 46 &\vvec 7{10} &\vvec 34  &\vvec 8{10}  & \vvec 56  &\vvec 78  &\vvec 22\\[10pt] 
\end{array}}
\quad \
\boxed{\begin{array}{ccccccccccccccccc}
~  \\[-10pt]
\vvec 12  &               &             &                    &                     &               &               &            &\vvec{\z2}{\z2}\\[14pt] 
\vvec 12  &               &             &                    &\vvec {\z1}{0}  &               &               &            &\vvec {\z2}{\z2}\\[14pt] 
\vvec 12  &               &\vvec 02 &                    &\vvec {\z1}{0}  &               & \vvec {\z3}{\z2}  &            &\vvec {\z2}{\z2}\\[14pt] 
\vvec 12  &\vvec14    &\vvec 02 &\vvec{\z1}2     &\vvec {\z1}{0} &\vvec{\z4}{\z2}    & \vvec {\z3}{\z2}  &\vvec{\z5}{\z4} &\vvec {\z2}{\z2}\\[10pt] 
\end{array}}
$$
\caption{Examples of Stern-Brocot arrays}
\label{fig:SBn}
\end{figure}  

\noindent
Spinor $\mathbf u = [ m\; n]^T$ determines curvature $curv = m^2+n^2 - 2$
via \eqref{eq:use}; 
for instance the last row of the array gives the following curvatures:  
$$
3,\quad 87, \quad 50, \quad 147, \quad 23, \quad 162, \quad 59,\quad 111, \quad 6
$$
(easy to locate in Figure \ref{fig:quarter2}).
Eventually, in the limit, one obtains all spinors in the range between the two initial ones.
Hence, the area of the corresponding fragment of the corona is 
$$
Area = \sum   \frac{\pi }{(m^2+n^2-2)^2} 
$$
where the sum is over all $(m,n)$ that appear in the Stern-Brocot produced by the above tree.
The sum of terms over a particular row $R_n$ gives the area of the circles up to the $n$-th depth.

To recover the full corona around the disk one may want to repeat the same process with 
one of the initial spinors changed to its negative, to account for the complementary segment 
of the corona (see the box on the right in Figure \ref{fig:SBn}). 
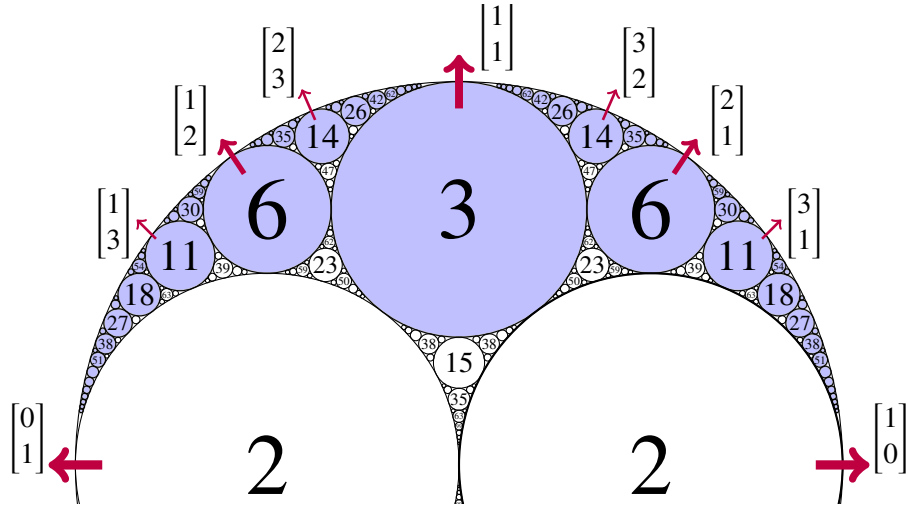
\begin{figure}
\centering
\begin{tikzpicture}[scale=5.1]
\clip (-1.25,-.1) rectangle (1.25,1.2);
\draw (0,0) circle (1);
%
\draw [thick]   (1/2 ,0/2) circle (1/2);
\draw [fill=blue, opacity=\sz, thick]   (0/3 ,2/3) circle (1/3);

\foreach \a/\b/\c/\d in {
3 / 4 /6,          
8 / 6 / 11,         
5 / 12/ 14,         
15/ 8 / 18,         
7 / 24 / 26,         
24/	10/	27,         
21/	20/	30,         
16/	30/	35,         
35/	12/	38,         
9/	40/	42,         
48/	14/	51,         
45/	28/	54,         
40/	42/	59,         
11/	60/	62,         
33/	56/	66,         
63/	16/	66,         
55/	48/	74,         
24/	70/	75,         
80/	18/	83,
13/	84/	86,         
77/	36/	86,         
39/	80/	90,         
65/	72/	98,         
99/	20/	102,
56/	90/	107,         
91/	60/	110,         
15/	112/	114,         
120/	22/	123,
117/	44/	126,         
32/	126/	131,         
112/	66/	131,         
105/	88/	138,         
17/	144/	146,
143/	24/	146,
96/	110/	147,         
51/	140/	150,         
85/	132/	158,         
119/	120/	170,         
72/	154/	171 ,        
168/	26/	171,
165/	52/	174,
160/	78/	179,
19/	180/	182,
57/	176/	186,
153/	104/	186,
95/	168/	194,
144/	130/	195,
195/	28/	198
}
\draw  [fill=blue, opacity=\sz, thick] (\a/\c,\b/\c) circle (1/\c)       (-\a/\c,\b/\c) circle (1/\c)
          (\a/\c,-\b/\c) circle (1/\c)       (-\a/\c,-\b/\c) circle (1/\c) 
;

\foreach \a/\b/\c/\s   in {
1 / 0 / 2 / 3
}
\draw (\a/\c,\b/\c) circle (1/\c)
          (-\a/\c,\b/\c) circle (1/\c)
node  [scale=\s]  at (\a/\c,\b/\c)   {$\c$}
node  [scale=\s]  at (-\a/\c,\b/\c) {$\c$};

\foreach \a/\b/\c/\s in {
0 / 2 / 3   /3,
0 /4 /15   /1,
0 / 6 / 35 / .6 
}
\draw (\a/\c,\b/\c) circle (1/\c)
          (\a/\c,-\b/\c) circle (1/\c)
        node [scale= \s]  at (\a/\c,\b/\c)  {$\c$}
        node [scale= \s]  at (\a/\c,-\b/\c) {$\c$};

\foreach \a/\b/\c in {
0 / 8/ 63,
0 /10 / 99
}
\draw (\a/\c,\b/\c) circle (1/\c)  
           (\a/\c,-\b/\c) circle (1/\c) 
node [scale= 20/\c]  at (\a/\c,\b/\c)  {$\c$}
node [scale= 20/\c]  at (\a/\c,-\b/\c) {$\c$};

\foreach \a/\b/\c in {
0 / 12 / 143,
0 / 14 / 195,
0 / 16 / 255,
0 / 18 / 323
}
\draw (\a/\c,\b/\c) circle (1/\c)  
           (\a/\c,-\b/\c) circle (1/\c);

\foreach \a/\b/\c/\s in {
3 / 4 /6/3,
8 / 6 / 11/1.5,
5 / 12/ 14/1.2
}
\draw (\a/\c,\b/\c) circle (1/\c)       (-\a/\c,\b/\c) circle (1/\c)
          (\a/\c,-\b/\c) circle (1/\c)       (-\a/\c,-\b/\c) circle (1/\c)
node [scale=\s]  at (\a/\c,\b/\c)    {$\c$}
node [scale=\s]  at (-\a/\c,\b/\c)   {$\c$}
node [scale=\s]  at (\a/\c,-\b/\c)   {$\c$}
node [scale=\s]  at (-\a/\c,-\b/\c)  {$\c$};

\foreach \a /  \b / \c   in {
15/ 8 / 18,
8 / 12 / 23,
7 / 24 / 26,
24/	10/	27,
21/	20/	30,
16/	30/	35,
3/	12/	38,
35/	12/	38,
24/	20/	39,
9/	40/	42,
16/	36/	47,
15/	24/	50,
48/	14/	51,
45/	28/	54,
24/	30/	59,
40/	42/	59,
11/	60/	62,
21/	36/	62,
48/	28/	63
}
\draw (\a/\c,\b/\c) circle (1/\c)          (-\a/\c,\b/\c) circle (1/\c)
          (\a/\c,-\b/\c) circle (1/\c)         (-\a/\c,-\b/\c) circle (1/\c)
node [scale=20/\c]  at (\a/\c,\b/\c)    {$\c$}
node [scale=20/\c]  at (-\a/\c,\b/\c)   {$\c$}
node [scale=20/\c]  at (\a/\c,-\b/\c)   {$\c$}
node [scale=20/\c]  at (-\a/\c,-\b/\c)  {$\c$};


\foreach \a /  \b / \c   in {
33/	56/	66,
63/	16/	66,
8/	24/	71,
55/	48/	74,
24/	70/	75,
48/	42/	83,
80/	18/	83,
13/	84/	86,
77/	36/	86,
24/	76/	87,
24/	40/	87,
39/	80/	90,
64/	60/	95,
80/	36/	95,
33/	72/	98,
63/	48/	98,
65/	72/	98,
3/	20/	102,
99/	20/	102,
8/	30/	107,
56/	90/	107,
5/	36/	110,
69/	60/	110,
91/	60/	110,
48/	92/	111,
48/	56/	111,
15/	112/	114,
15/	40/	114,
40/	72/	119,
39/	96/	122,
120/	22/	123,
117/	44/	126,
32/	126/	131,
48/	102/	131,
112/	66/	131,
35/	60/	134,
99/	60/	134,
120/	44/	135,
105/	88/	138,
32/	132/	143,
80/	72/	143,
17/	144/	146,
143/	24/	146,
48/	70/	147,
96/	110/	147,
45/	76/	150,
51/	140/	150,
120/	66/	155,
85/	132/	158,
63/	80/	162,
24/	60/	167,
56/	120/	167,
136/	84/	167,
55/	96/	170,
119/	120/	170,
72/	154/	171,
168/	26/	171,
69/	92/	174,
165/	52/	174,
64/	102/	179,
80/	90/	179,
160/	78/	179,
19/	180/	182,
51/	156/	182,
141/	84/	182,
120/	88/	183,
168/	52/	183,
57/	176/	186,
153/	104/	186,
48/	84/	191,
128/	132/	191,
65/	120/	194,
95/	168/	194,
129/	120/	194,
144/	130/	195,
3/	28/	198,
195/	28/	198,
8/	42/	203
}
\draw (\a/\c,\b/\c) circle (1/\c)          (-\a/\c,\b/\c) circle (1/\c)
          (\a/\c,-\b/\c) circle (1/\c)         (-\a/\c,-\b/\c) circle (1/\c);

\draw [->, line width=3.5pt, purple]  (.93,0) -- (1.07,0);    
\draw [->, line width=3.5pt, purple]  (-.93,0) -- (-1.07,0);   

\draw [->, line width=3pt, purple] (0,.93) -- (0,1.07);     

\draw [->, line width=2pt, purple] (.56,.76) -- (0.62,0.85);
\draw [->, line width=2pt, purple] (-.56,.76) -- (-0.62,0.85);

\draw [->, line width=1pt, purple] (.79,.59) -- (0.84,0.64);
\draw [->, line width=1pt, purple] (-.79,.59) -- (-0.84,0.64);

\draw [->, line width=1pt, purple] (.375,.9) -- (0.41,0.98);
\draw [->, line width=1pt, purple] (-.375,.9) -- (-0.41,0.98);

\node at (1.125, 0.07) {$\begin{bmatrix} 1\\ 0 \end{bmatrix}$};
\node at (-1.125, 0.07) {$\begin{bmatrix} 0\\ 1 \end{bmatrix}$};
\node at (0.1, 1.12) {$\begin{bmatrix} 1\\ 1 \end{bmatrix}$};
\node at (0.7, .9) {$\begin{bmatrix} 2\\ 1 \end{bmatrix}$};
\node at (-0.7, .9) {$\begin{bmatrix} 1\\ 2 \end{bmatrix}$};

\node at (0.9, .63) {$\begin{bmatrix} 3\\ 1 \end{bmatrix}$};
\node at (-0.9, .63) {$\begin{bmatrix} 1\\ 3 \end{bmatrix}$};

\node at (0.47, 1.05) {$\begin{bmatrix} 3\\ 2 \end{bmatrix}$};
\node at (-0.47, 1.05) {$\begin{bmatrix} 2\\ 3 \end{bmatrix}$};

\end{tikzpicture}
\caption{Apollonian Window -- right op corner fragment  the major corona}
\label{fig:windowspin}
\end{figure}

The sum would have to be over a set of vectors
reproduced case-by-case.
This is sufficient for estimation but it lacks an elegant algebraic form.
But there is a remedy, as presented in the following.

~

\noindent
{\bf B. Universal case.}  
The upper half of the corona in Apollonian Window  
around the great circle of curvature $(-1)$, shown in Figure \ref{fig:windowspin}, 
consists of spinors forming the following Stern-Brocot array (the first five rows are shown):
\begin{figure}[H]
{\small
$$
\arraycolsep=2.4pt
\begin{array}{ccccccccccccccccc}
\vvec 01  &&               &&             &&                &&              &&               &&               &&            &&\vvec 10\\[12pt] 
\vvec 01  &&               &&             &&                &&\vvec 11  &&               &&               &&            &&\vvec 10\\[12pt] 
\vvec 01  &&               &&\vvec 12 &&                &&\vvec 11  &&               && \vvec 21  &&            &&\vvec 10\\[12pt] 
\vvec 01  &&\vvec13    &&\vvec 12 &&\vvec23     &&\vvec 11  &&\vvec32    && \vvec 21  &&\vvec31 &&\vvec 10\\[12pt] 
\vvec 01  &\vvec 14&\vvec13    &\vvec 25&\vvec 12 &\vvec 35&\vvec23   & \vvec 34 &\vvec 11 & \vvec 43&\vvec32  
& \vvec 54 & \vvec 21 &\vvec 42 &\vvec31 & \vvec 41  &\vvec 10 
\end{array}
$$
}
\caption{Universal Stern-Brocot array}
\label{fig:universal}
\end{figure}
This arrangement proves to be {\bf universal} for our purpose.
It has two essential properties:  
(1) the coefficients $k$ and $n$ of every entry $[k\, n]^T$ are co-primes, and   
(2) every coprime pair appears in the tree.
This follows from the fact that if 
we interpret the vector as a fraction, $\vvec kn \mapsto k/n$, 
the above tree becomes the well-known 
classical Stern-Brocot tree of positive rational numbers, 
which is known to lists them in the form of reduced fractions, each exactly once
(one must ignore the repetitions down along the columns).

The vectors of Stern-Brocot tree of spinors 
may be visualized as so-called Euclid'd orchard (see Fig \ref{fig:Euclid}, left).
To account also for the bottom part of the Apollonian Window,
we must extend the set to include spinors produced by 
the pair 
$$
\begin{bmatrix} 0 \\ 1  \end{bmatrix},   
\quad
\begin{bmatrix} -1 \\\phantom{-}0  \end{bmatrix} \,.
$$
This would graphically correspond
to the ``extended orchard'' shown in Figure \ref{fig:Euclid}, center.
Note however that one entry is repeated, since 
$[-1,0]^T\sim[1,0]^T$ and both correspond to the same disk.
Instead of removing one copy, extend the orchard 
to all pairs of co-primes in the lattice $\mathbb Z\times \mathbb Z$,
as shown in Figure \ref{fig:Euclid}, right.
Now every entry appears exactly twice. We have take half of the sum,
but we gain more elegant formulation:

\begin{figure}[t]
\centering
\includegraphics[scale=.8]{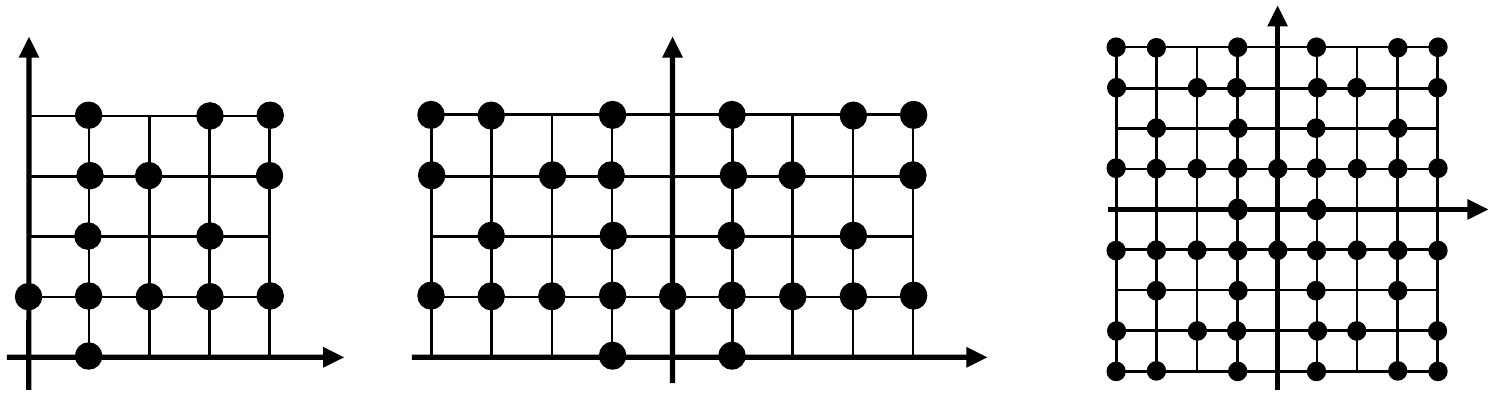}
\caption{\small Euclidean orchard and its extensions.}
\label{fig:Euclid}
\end{figure}


\begin{proposition}
\label{thm:A}
The total area of disks tangent to the great circle in the Apollonian window is 
$$
A \ = \ \frac{\pi}{2}   \sum_{\substack{m,n\,\in \mathbb Z \\ \gcd(m,n)=1}} \ 
                                             \frac{1}{ \big( m^2 + n^2 +1\big)^2}
$$
where $m$ and $n$ run over all pairs coprime integers (including negative).
\end{proposition}

Denote  $\mathbb  Z^2_{\Euc}$ the subset of the integral lattice $\mathbb Z\times\mathbb Z$ 
restricted to vectors with co-prime coefficients:
\begin{equation}
\renewcommand*{\arraystretch}{.7}
\setlength\arraycolsep{2pt}
\mathbb Z^2_{\Euc} \ = \ \left\{ \; \vvec{m}{n}  \in \mathbb Z \times \mathbb Z \;\big|\; \gcd(m,n) = 1\,   \right\} 
\end{equation}
(it can be called the ``extended Euclid's orchard''). 
The above formula may now be rewritten as 
$$
A =  \frac{1}{2} \sum_{\mathbf v \in \mathbb Z^2_{\Euc}} \ \frac{1} {\big( \; \|\mathbf v \|^2 +1\big)^2}
$$

~

\noindent
{\bf C.  General case.}
Now, to get a formula for an arbitrary corona,
notice that the Stern-Brocot tree determined
by vectors
$$
\mathbf u=\begin{bmatrix} a \\ b  \end{bmatrix}, \ \  
\mathbf v=\begin{bmatrix} c \\ d  \end{bmatrix}, \quad 
$$
may be obtained from the universal Stern-Brocot tree (Figure \ref{fig:universal})
by replacing its every entry
by a vector transformed by matrix
\begin{equation}
\label{eq:MMM}
M \ = \  [\mathbf u|\mathbf v] \  \equiv \ \begin{bmatrix} a & c  \\ b & d \end{bmatrix}
\ \in \ {\rm End} (n, \mathbb R) \ \cong \  {\rm Mat}_{2\times 2}(\mathbb R)
\end{equation}
For instance, the Stern-Brocot array in Figure \ref{fig:SBn} may be viewed as the product of the map
$$
\mathbb Z^2_o \ \ni \  \begin{bmatrix} m \\ n  \end{bmatrix}
\quad  \mapsto \quad  
\begin{bmatrix} 1&2 \\ 2&2  \end{bmatrix}
\begin{bmatrix} m \\ n  \end{bmatrix}
\ \in \  \mathbb R^2
$$ 
%
%
%
%
Hence our main result:

~

\noindent
{\bf The Main Theorem \  [Same as Thm. \ref{thm:f}]:}  The area of the corona of the disk $B$ 
with two adjacent spinors forming matrix $M$ (see \eqref{eq:MMM}) is 
\begin{equation}
\label{eq:main2}
\hbox{\rm Entire corona:}\qquad 
         K(M) \ = \ \frac{1}{2} \ \sum_{\mathbf f \in \mathbb Z^2_{\Euc}} \ 
                        \big( \langle\; \mathbf f^T \,|\, M ^TM \,|\,  \mathbf f \;\rangle  - B\big)^{-2}
\end{equation}

~

To account for a fragment of corona between the circles related to spinors $\mathbf u$ and $\mathbf v$,
we use the same formula except with the summation going over one of these sets
$$
\begin{aligned}
(a)\quad & \mathbb N^2_{\sf o} \ =: \ \{\, (m,n)\in \mathbb N\times \mathbb N \;|\; \gcd(m,n)=1 \,\} \\
(b)\quad & \dot{\mathbb N}^2_{\sf o} \ =: \ \{\, (m,n)\in \mathbb N\times \mathbb N \;|\; \gcd(m,n)=1 \,\}
\end{aligned}
$$
where $\mathbb N = \{1,2,...\}$ and $\dot{\mathbb N}_{\sf o}  = \{0,1,2,...\}$.
Case (a)  captures disks strictly between the end-disks, while (b) includes the end disks
(see Figure \ref{fig:Euclid}, left).

~

The formula for the total corona area $A$ involves a choice two spinors,
yet its value is is invariant with respect to this choice.
Indeed,
the group ${\rm SL(}2,\mathbb Z)$ acts on $\mathbb Z^2_o$ as bijection
and permutes spinors, preserving their mutual relations of neighborhood.
Thus the matrix $M$ is transformed accordingly:
$$
M = [\mathbf u | \mathbf v] 
\ \mapsto \ 
[\mathbf u | \mathbf v]g =   Mg
$$

\begin{proposition}(Invariance).  
\label{thm:invariance}
The area function of the corona is invariant under action of the modular group, that is:
$$
K(M) = K(gM)    \qquad \hbox{\rm for any} \  g\in {\rm SL}(2,\mathbb Z)\,.
$$
\end{proposition}

\noindent
{\bf Proof:} Denote the summand in \eqref{eq:main} by $A(M,\mathbf f, B)$.
Then 

\def\dsum{\displaystyle\sum}

$$
\begin{array}{rl}
 \dsum_{\mathbf f \,\in\, \mathbb Z^2_{\Euc}}  A\left[ (Mg)\mathbf f\right]
\ &= \ 
\dsum_{\mathbf f \,\in\, \mathbb Z^2_{\Euc}}  A\left[M(g\mathbf f)\right]
\ = \ 
\dsum_{g\mathbf f \,\in\, \mathbb Z^2_{\Euc}}  A\left[M\mathbf f\right]  \\[19pt]
\ &= 
\dsum_{\mathbf f \,\in\, g^{-1}(\mathbb Z^2_{\Euc})}  A\left[M\mathbf f\right]
\ = \ 
\dsum_{\mathbf f \,\in\, \mathbb Z^2_{\Euc}}  A\left[M\mathbf f\right]
\end{array}
$$
\qed

~~

The main theorem is equivalent to the following statement:

\begin{proposition}
Let $\mathbf u$ and $\mathbf v$ be two neighboring spinors in a corona around a disk of curvature $B$ 
in an Apollonian disk packing.  Maintain the previous notation:
$$
\mathbf u=\begin{bmatrix} a \\ b  \end{bmatrix}, \ \  
\mathbf v=\begin{bmatrix} c \\ d  \end{bmatrix}, \quad 
\quad
M \ = \  [\mathbf u|\mathbf v] \  \equiv \ \begin{bmatrix} a & c  \\ b & d \end{bmatrix},
\quad
\mathbf f = \begin{bmatrix} m \\ n  \end{bmatrix}
$$ 
Then the quadratic polynomial
$$
P(m,n) \ = \ \mathbf f^T(M^TM)\;\mathbf f - B
$$
reproduces all curvatures in the corona of $B$ with the integers $m,n\in \mathbf Z^2$ running of all pairs of integer co-primes,
including $(\pm 1,0)$ and $(0,\pm 1)$.
Each curvature appears in the range of the polynomial twice due to $P(\mathbf f) = P(-\mathbf f)$.
\end{proposition}
In terms of the entries of the matrices, the polynomial is 
$$
P(m,n) \ = \  (a^2\!+\!c^2)\, m^2  \ +\  (c^2\!+\!d^2)\, n^2 + (ab\!+\!cd)\, mn \ - \ B
$$
The invariance of the range of the above polynomial 
with respect to the group $\SL(2,Z$ holds 
for the same argument as in Proposition \ref{thm:invariance}.
\newpage

\section{Examples and special cases}

\begin{figure}[H]
\begin{tikzpicture}[scale=4.5]  

\clip (-.2,-.2) rectangle (1.5,2.15);

\foreach \a/\b/\c/\d in {
3/  4/ 4 ,
8/    6/   9 ,
24/  20/  25,
24/  10/  25,
15/   8/   16, 
35/	12/	36,
48/	42/	49,
48/	28/	49,
48/	14/	49,
63/	48/	64,	
63/	16/	64,	
80/	72/	81,	
80/	36/	81,	
80/	18/	81,	
99/	60/	100,	
99/	20/	100,	
120/	110/	121,	
120/	88/	121,	
120/	66/	121,	
120/	44/	121,	
120/	22/	121,
143/	120/	144,	
143/	24/	144,
168/	156/	169,	
168/	130/	169,
168/	104/	169,	
168/	78/	169,	
168/	52/	169,	
168/	26/	169,	
189/	148/	196,	
195/	140/	196,	
195/	84/	196,	
195/	28/	196	
}
\draw  [fill=blue, opacity=\sz] (\a/\c,\b/\c) circle (1/\c)    
                                                    (\a/\c,2-\b/\c) circle (1/\c)   
;
\draw (1,-1) -- (1,3);
\draw (-1,-1) -- (-1,3);
\draw (0,0) circle (1);
\draw [thick] (0,0) circle (1)
          node  [scale=5]  at (0,0) {$1$};
\draw (0,2) circle (1)
          node  [scale=4]  at (0,1.8) {$1$};
\foreach \a/\b/\c/\s in {
3/  4/ 4/    3,
5/ 12/12/  1, 
7/ 24/24/ .8,
9/ 40/40/ .6
}
\draw (\a/\c,\b/\c) circle (1/\c)    (\a/\c,-\b/\c) circle (1/\c)
          (-\a/\c,\b/\c) circle (1/\c)    (-\a/\c,-\b/\c) circle (1/\c)
        node [scale= \s]  at (\a/\c,\b/\c)  {$\c$}
        node [scale= \s]  at (\a/\c,-\b/\c) {$\c$}
        node [scale= \s]  at (-\a/\c,\b/\c) {$\c$}
        node [scale= \s]  at (-\a/\c,-\b/\c) {$\c$};
\foreach \a/\b/\c in {
11/	60	/60,
13/	84/	84,
15/	112/	112
}
\draw (\a/\c,\b/\c) circle (1/\c)  
           (-\a/\c,\b/\c) circle (1/\c) 
node [scale= 20/\c]  at (\a/\c,\b/\c)  {$\c$}
node [scale= 20/\c]  at (\a/\c,-\b/\c) {$\}$}
node [scale= 20/\c]  at (-\a/\c,\b/\c) {$\c$}
node [scale= 20/\c]  at (-\a/\c,-\b/\c) {$\c$};

\foreach \a/\b/\c in {
17/	144/	144,
19/	180/	180,
21/	220/	220,
23/	264/	264,
25/	312/	312
}
\draw (\a/\c,\b/\c) circle (1/\c)   (\a/\c, -\b/\c) circle (1/\c)
          (-\a/\c, \b/\c) circle (1/\c)   (-\a/\c, -\b/\c) circle (1/\c);

\foreach \a/\b/\c/\d in {
8/    6/   9     / 2,
15/   8/   16  / 1
}
\draw (\a/\c,\b/\c) circle (1/\c)       (-\a/\c,\b/\c) circle (1/\c)
          (\a/\c,-\b/\c) circle (1/\c)       (-\a/\c,-\b/\c) circle (1/\c)
          (\a/\c,2-\b/\c) circle (1/\c)       (-\a/\c,2-\b/\c) circle (1/\c)
node [scale=\d]  at (\a/\c,\b/\c)    {$\c$}
node [scale=\d]  at (-\a/\c,\b/\c)   {$\c$}
node [scale=\d]  at (\a/\c,-\b/\c)   {$\c$}
node [scale=\d]  at (-\a/\c,-\b/\c)  {$\c$}
node [scale=\d]  at (\a/\c,2-\b/\c)    {$\c$}
node [scale=\d]  at (-\a/\c,2-\b/\c)   {$\c$};

\foreach \a /  \b / \c   in {
24/20/25,   24/10/25,    21/20/28,
16/30/33,   35/12/36,    48/42/49,
48/28/49,   48/14/49,    45/28/52
}
\draw (\a/\c,\b/\c) circle (1/\c)          (-\a/\c,\b/\c) circle (1/\c)
          (\a/\c,-\b/\c) circle (1/\c)         (-\a/\c,-\b/\c) circle (1/\c)
          (\a/\c,2-\b/\c) circle (1/\c)          (-\a/\c,2-\b/\c) circle (1/\c)
node [scale=18/\c]  at (\a/\c,\b/\c)    {$\c$}
node [scale=18/\c]  at (-\a/\c,\b/\c)   {$\c$}
node [scale=18/\c]  at (\a/\c,-\b/\c)   {$\c$}
node [scale=18/\c]  at (-\a/\c,-\b/\c)  {$\c$}
node [scale=18/\c]  at (\a/\c,2-\b/\c)    {$\c$}
node [scale=18/\c]  at (-\a/\c,2-\b/\c)   {$\c$};


\foreach \a /  \b / \c   in {
40/	42/	57,	33/	56/	64,	63/	48/	64,	63/	16/	64,	55/	48/	72,	24/	70/	73,
69/	60/	76,	80/	72/	81,	64/	60/	81,	80/	36/	81,	80/	18/	81,	77/	36/	84,
39/	80/	88,	65/	72/	96,	48/	92/	97,	99/	60/	100,	99/	20/	100,	56/	90/	105,
91/	60/	108,	120/	110/	121,	120/	88/	121,	120/	66/	121,	120/	44/	121,	120/	22/	121,
117/	44/	124,	32/	126/	129,	112/	66/	129,	105/	88/	136,	143/	120/	144,	143/	24/	144,
96/	110/	145,	51/	140/	148,	141/	84/	148,	136/	120/	153,	136/	84/	153,	149/	132/	156,
85/	132/	156,	129/	120/	160,	119/	120/	168,	168/	156/	169,	72/	154/	169,	168/	130/	169,
168/	104/	169,	168/	78/	169,	168/	52/	169,	168/	26/	169,	165/	52/	172,	128/	132/	177,
160/	78/	177,	57/	176/	184,	153/	104/	184,	95/	168/	192,	96/	186/	193,	144/	130/	193,
189/	148/	196,	195/	140/	196,	195/	84/	196,	195/	28/	196,	40/	198/	201,	104/	180/	201,
184/	162/	201
}
\draw (\a/\c,\b/\c) circle (1/\c)          (-\a/\c,\b/\c) circle (1/\c)
          (\a/\c,-\b/\c) circle (1/\c)         (-\a/\c,-\b/\c) circle (1/\c)
          (\a/\c,2-\b/\c) circle (1/\c)          (-\a/\c,2-\b/\c) circle (1/\c);

\draw [->, line width=3.2pt, purple]  (.83,0) -- (1.17,0);
\node at (1.37, 0) [scale=1.2] {$\begin{bmatrix} 1 \  0\, \end{bmatrix}$};

\draw [->, line width=3.2pt, purple]  (.83,2) -- (1.17,2);
\node at (1.37, 2) [scale=1.2] {$\begin{bmatrix} 1 \ 0\, \end{bmatrix}$};

\draw [->, line width=3pt, purple] (.87,1) -- (1.14,1);
\node at (1.32, 1) [scale=1.1] {$\begin{bmatrix} 2\ 0\, \end{bmatrix}$};

\draw [->, line width=2.5pt, purple]  (.93,.667) -- (1.08,.667);
\node at (1.23, .667) {$\begin{bmatrix} 3\  0\, \end{bmatrix}$};

\draw [->, line width=2.5pt, purple]  (.93, 1.333) -- (1.08, 1.333);
\node at (1.23, 1.333) {$\begin{bmatrix} 3 \  0\, \end{bmatrix}$};

\draw [->, line width=2pt, purple]  (.98, 1.2) -- (1.05, 1.2);
\node at (1.18, 1.2) [scale=.8]{$\begin{bmatrix} 5 \  0\, \end{bmatrix}$};

\draw [->, line width=2pt, purple]  (.98, 0.8) -- (1.05, 0.8);
\node at (1.18, 0.8)  [scale=.8] {$\begin{bmatrix} 5 \ 0\, \end{bmatrix}$};

\draw [->, line width=2pt, purple]  (.98, 0.5) -- (1.05, 0.5);
\node at (1.18, 0.5) [scale=.8] {$\begin{bmatrix} 4 \ 0\,  \end{bmatrix}$};

\draw [->, line width=2pt, purple]  (.98, 1.5) -- (1.05, 1.5);
\node at (1.18, 1.5) [scale=.8] {$\begin{bmatrix} 4 \ 0\, \end{bmatrix}$};

\end{tikzpicture}
\caption{Ford circles as a corona in the Apollonian Belt}
\label{fig:belt}
\end{figure}
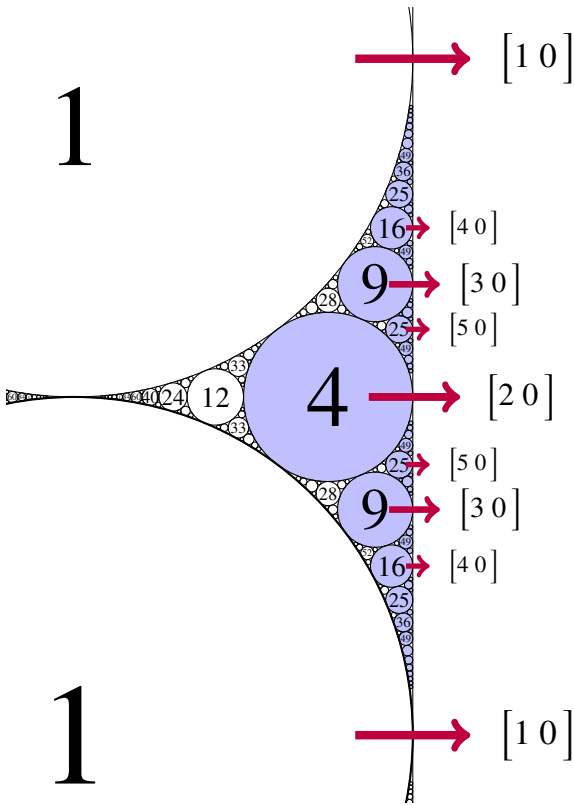

\noindent
{\bf Example 1: Ford circles}.
Set  the Apollonian Belt vertically, as in Figure \ref{fig:belt}.  
The spinors at the two unit circles are 
$$
\mathbf v = \begin{bmatrix} 1\\ 0\end{bmatrix}
\qquad\hbox{and}\qquad
\mathbf w = \begin{bmatrix} 1\\ 0\end{bmatrix}
$$
The general formula \eqref{eq:main} becomes
$$
A \ = \ \pi \sum_{\gcd(m,k)=1} \ \frac{1}{ (m+k)^4}  
   \ = \ \pi \sum_{n=1}^\infty \ \frac{\varphi(n)}{ n^4} 
$$
where $\varphi(n)$ denotes the Euler's totient function 
(the number of positive coprimes with $n$ not exceeding $n$).
The last equation  is implied by a simple -- easy to show -- property true for $n>1$:
$$
\varphi(n) = |\,\{(i,j)\in \mathbb N\times \mathbb N \;\big|\; i+j=n,\, \gcd(i,j) = 1\}\,|
$$
The result is an example of Dirichlet series, which is known to converge to a ratio of Riemann zeta functions \cite{HW}.
In general:
$$
\sum_{n=1}^\infty \ \frac{\varphi(n)}{ n^s} = \frac{\zeta(s-1)}{\zeta(s)}
$$
%


\noindent
Hence, in our case:
$$
A_{\rm Ford}/\pi  \ = \ 1 + \sum_{n=1}^\infty \ \frac{\varphi(n)}{ n^4}   \ = \ 1+ \frac{\zeta(3)}{\zeta(4)} 
$$

\begin{figure}
\begin{tikzpicture}[scale=6]  
\clip (-.05,-.1) rectangle (1.25,1.2);
\draw (0,0) circle (1);

\draw [fill=blue, opacity=\sz, thick]   (1/2 ,0/2) circle (1/2);
\draw [fill=blue, opacity=\sz, thick]   (0/3 ,2/3) circle (1/3);

\foreach \a/\b/\c/\d in {
3 / 4 /6,          
8 / 6 / 11,         
5 / 12/ 14,         
15/ 8 / 18,         
7 / 24 / 26,         
24/	10/	27,         
21/	20/	30,         
16/	30/	35,         
35/	12/	38,         
9/	40/	42,         
48/	14/	51,         
45/	28/	54,         
40/	42/	59,         
11/	60/	62,         
33/	56/	66,         
63/	16/	66,         
55/	48/	74,         
24/	70/	75,         
80/	18/	83,
13/	84/	86,         
77/	36/	86,         
39/	80/	90,         
65/	72/	98,         
99/	20/	102,
56/	90/	107,         
91/	60/	110,         
15/	112/	114,         
120/	22/	123,
117/	44/	126,         
32/	126/	131,         
112/	66/	131,         
105/	88/	138,         
17/	144/	146,
143/	24/	146,
96/	110/	147,         
51/	140/	150,         
85/	132/	158,         
119/	120/	170,         
72/	154/	171 ,        
168/	26/	171,
165/	52/	174,
160/	78/	179,
19/	180/	182,
57/	176/	186,
153/	104/	186,
95/	168/	194,
144/	130/	195,
195/	28/	198
}
\draw  [fill=blue, opacity=\sz, thick] (\a/\c,\b/\c) circle (1/\c)       (-\a/\c,\b/\c) circle (1/\c)
          (\a/\c,-\b/\c) circle (1/\c)       (-\a/\c,-\b/\c) circle (1/\c) 
;

\foreach \a/\b/\c/\s   in {
1 / 0 / 2 / 3
}
\draw (\a/\c,\b/\c) circle (1/\c)
          (-\a/\c,\b/\c) circle (1/\c)
node  [scale=\s]  at (\a/\c,\b/\c)   {$\c$}
node  [scale=\s]  at (-\a/\c,\b/\c) {$\c$};

\foreach \a/\b/\c/\s in {
0 / 2 / 3   /3,
0 /4 /15   /1,
0 / 6 / 35 / .6 
}
\draw (\a/\c,\b/\c) circle (1/\c)
          (\a/\c,-\b/\c) circle (1/\c)
        node [scale= \s]  at (\a/\c,\b/\c)  {$\c$}
        node [scale= \s]  at (\a/\c,-\b/\c) {$\c$};

\foreach \a/\b/\c in {
0 / 8/ 63,
0 /10 / 99
}
\draw (\a/\c,\b/\c) circle (1/\c)  
           (\a/\c,-\b/\c) circle (1/\c) 
node [scale= 20/\c]  at (\a/\c,\b/\c)  {$\c$}
node [scale= 20/\c]  at (\a/\c,-\b/\c) {$\c$};

\foreach \a/\b/\c in {
0 / 12 / 143,
0 / 14 / 195,
0 / 16 / 255,
0 / 18 / 323
}
\draw (\a/\c,\b/\c) circle (1/\c)  
           (\a/\c,-\b/\c) circle (1/\c);

\foreach \a/\b/\c/\s in {
3 / 4 /6/3,
8 / 6 / 11/1.5,
5 / 12/ 14/1.2
}
\draw (\a/\c,\b/\c) circle (1/\c)       (-\a/\c,\b/\c) circle (1/\c)
          (\a/\c,-\b/\c) circle (1/\c)       (-\a/\c,-\b/\c) circle (1/\c)
node [scale=\s]  at (\a/\c,\b/\c)    {$\c$}
node [scale=\s]  at (-\a/\c,\b/\c)   {$\c$}
node [scale=\s]  at (\a/\c,-\b/\c)   {$\c$}
node [scale=\s]  at (-\a/\c,-\b/\c)  {$\c$};

\foreach \a /  \b / \c   in {
15/ 8 / 18,
8 / 12 / 23,
7 / 24 / 26,
24/	10/	27,
21/	20/	30,
16/	30/	35,
3/	12/	38,
35/	12/	38,
24/	20/	39,
9/	40/	42,
16/	36/	47,
15/	24/	50,
48/	14/	51,
45/	28/	54,
24/	30/	59,
40/	42/	59,
11/	60/	62,
21/	36/	62,
48/	28/	63
}
\draw (\a/\c,\b/\c) circle (1/\c)          (-\a/\c,\b/\c) circle (1/\c)
          (\a/\c,-\b/\c) circle (1/\c)         (-\a/\c,-\b/\c) circle (1/\c)
node [scale=20/\c]  at (\a/\c,\b/\c)    {$\c$}
node [scale=20/\c]  at (-\a/\c,\b/\c)   {$\c$}
node [scale=20/\c]  at (\a/\c,-\b/\c)   {$\c$}
node [scale=20/\c]  at (-\a/\c,-\b/\c)  {$\c$};


\foreach \a /  \b / \c   in {
33/	56/	66,
63/	16/	66,
8/	24/	71,
55/	48/	74,
24/	70/	75,
48/	42/	83,
80/	18/	83,
13/	84/	86,
77/	36/	86,
24/	76/	87,
24/	40/	87,
39/	80/	90,
64/	60/	95,
80/	36/	95,
33/	72/	98,
63/	48/	98,
65/	72/	98,
3/	20/	102,
99/	20/	102,
8/	30/	107,
56/	90/	107,
5/	36/	110,
69/	60/	110,
91/	60/	110,
48/	92/	111,
48/	56/	111,
15/	112/	114,
15/	40/	114,
40/	72/	119,
39/	96/	122,
120/	22/	123,
117/	44/	126,
32/	126/	131,
48/	102/	131,
112/	66/	131,
35/	60/	134,
99/	60/	134,
120/	44/	135,
105/	88/	138,
32/	132/	143,
80/	72/	143,
17/	144/	146,
143/	24/	146,
48/	70/	147,
96/	110/	147,
45/	76/	150,
51/	140/	150,
120/	66/	155,
85/	132/	158,
63/	80/	162,
24/	60/	167,
56/	120/	167,
136/	84/	167,
55/	96/	170,
119/	120/	170,
72/	154/	171,
168/	26/	171,
69/	92/	174,
165/	52/	174,
64/	102/	179,
80/	90/	179,
160/	78/	179,
19/	180/	182,
51/	156/	182,
141/	84/	182,
120/	88/	183,
168/	52/	183,
57/	176/	186,
153/	104/	186,
48/	84/	191,
128/	132/	191,
65/	120/	194,
95/	168/	194,
129/	120/	194,
144/	130/	195,
3/	28/	198,
195/	28/	198,
8/	42/	203
}
\draw (\a/\c,\b/\c) circle (1/\c)          (-\a/\c,\b/\c) circle (1/\c)
          (\a/\c,-\b/\c) circle (1/\c)         (-\a/\c,-\b/\c) circle (1/\c);

\draw [->, line width=3pt, purple]  (.93,0) -- (1.07,0);
\draw [->, line width=3pt, purple] (0,.93) -- (0,1.07);

\node at (1.125, 0.07) {$\begin{bmatrix} 1\\ 0 \end{bmatrix}$};
\node at (0.07, 1.1) {$\begin{bmatrix} 1\\ 1 \end{bmatrix}$};
\end{tikzpicture}
\caption{Apollonian Window -- right op corner fragment  the major corona}
\label{fig:quarter}
\end{figure}
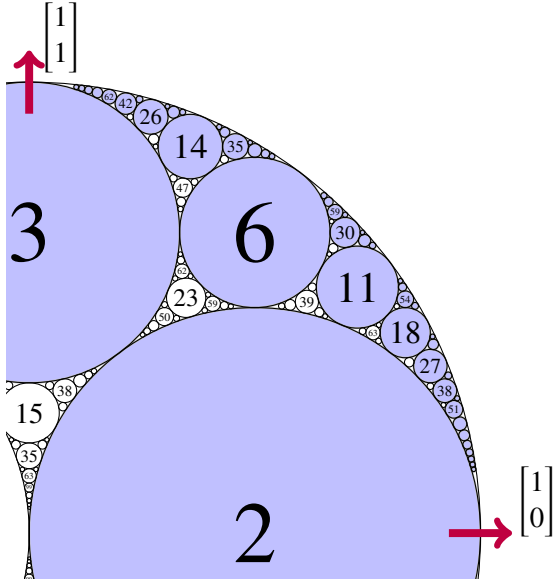

\noindent
{\bf Example 2:  Apollonian Window -- a quarter.}
Consider the upper right quarter of the corona of the great circle in the Apollonian Window,
Fig. \ref{fig:quarter}.
The spinors at the two big circle are 
$$
\mathbf v = \begin{bmatrix} 1\\ 1\end{bmatrix}\,,
\qquad
\mathbf w = \begin{bmatrix} 0\\ 1\end{bmatrix}\,,
\qquad
B=-1
$$
The general formula \eqref{eq:main2} implies
$$
A_{\mathbf v,\mathbf w} \ = \ \sum_{n=1}^\infty \sum_{\substack{k\leq n \\ (n,k)=1}} \ \frac{\pi}{\big( n^2+k^2+1\big)^2}
$$


\noindent
{\bf Example 3:  Apollonian Window -- around $b=2$.}
Consider the corona around the right disc of curvature $B=2$
(yellow disk in Fig. \ref{fig:quarter2}).
Following \eqref{eq:main}, the area of the segment of disks between the spinors chosen below is   
$$
\mathbf v = \begin{bmatrix} 1\\ 0\end{bmatrix}
\ \hbox{and}\ 
\mathbf w = \begin{bmatrix} 0\\ 2\end{bmatrix}
\qquad
\Rightarrow
\qquad
A_{\mathbf v,\mathbf w} \ = \  \sum_{n=1}^\infty \sum_{\substack{ k\leq n \\  (n,k)=1}}  \ \frac{1}{\big( n^2+4k^2-2\big)^2}
$$
Similarly, the fragment between the circles of curvature $3$ and $(-1)$ (the one that includes the chain 6, 11, 18 etc)
can be expressed by spinors 
$$
\mathbf v = \begin{bmatrix} 1\\ 0\end{bmatrix}
\qquad
\mathbf w = \begin{bmatrix} 1\\ 2\end{bmatrix}
\qquad
\qquad
\Rightarrow
\qquad
A_{\mathbf v,\mathbf w} \ = \ \sum_{\substack {k,n\in \mathbb N \\  (n,k)=1}}\ \frac{1}{\big( n^2+2nk+ 5k^2-2\big)^2}
$$
For the total area, 
starting with spinors $[0\;1]^T$ and $[2\;0]^T$, we get
$$
A \ = \ \frac{1}{2}  \sum_{\substack {k,n\in \mathbb N \\  \gcd(n,k)=1}}\ \frac{1}{\big( n^2+ 4k^2-2\big)^2}
$$

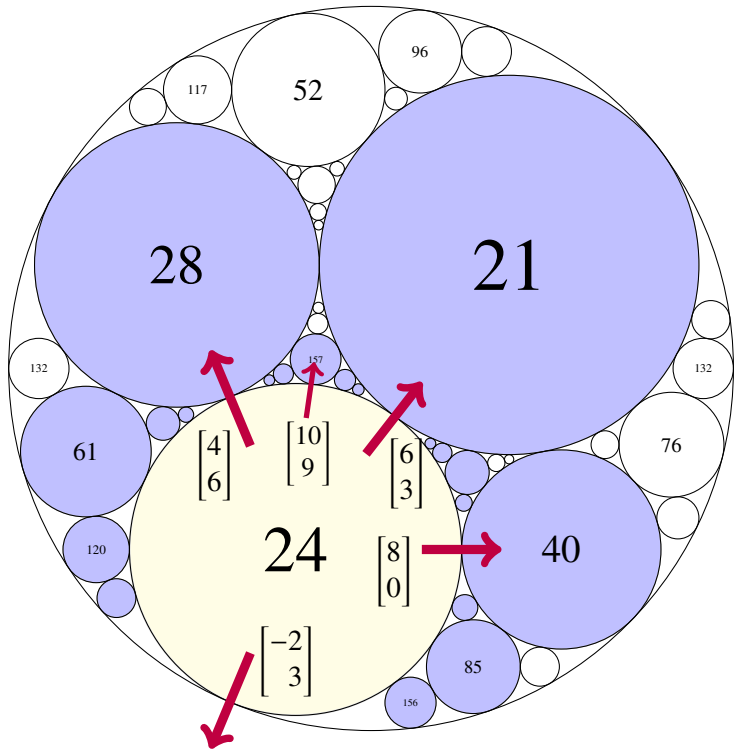
\begin{figure}[H]
\centering
\begin{tikzpicture}[scale=53]

\draw [fill=yellow, opacity=.1, thick] (17/24,12/24) circle (1/24);
%
\draw (8/11,6/11) circle (1/11);
\foreach \a /  \b / \c   in {
16/	12/	21,
19/	16/	28,
31/	20/	40,
40/	32/	61,
64/	40/	85,
79/	60/	120,
115	/72/	156,
112/	86/	157,
136/	94/	181,
136/	100/	205,
160/	126/	237,
232/	150/	309,
271/	204/	376,
280/	216/	397,
307/	216/	412,
352/	240/	469,
352/	276/	517,
496/	370/	685,
544/	386/	733,
520/	402/	741
}
\draw [fill=blue, opacity=\sz] (\a/\c,\b/\c) circle (1/\c);

\foreach \a /  \b / \c   in {
16/	12/	21,
17/	12/	24,
19/	16/	28,
31/	20/	40,
37/	32/	52,
40/	32/	61,
61/	40/	76,
64/	40/	85,
71/	60/	96,
80/	72/	117,
79/	60/	120,
85/	72/	132,
107 /72/	132,
115	/72/	156,
112/	86/	157
}
\draw (\a/\c,\b/\c) circle (1/\c)         
node [scale=54/\c]  at (\a/\c,\b/\c)    {$\c$}
;


\foreach \a /  \b / \c   in {
121/	100/	160,
136/	94/	181,
152/	96/	189,
157/	96/	204,
136/	100/	205,
169/	116/	208,
152/	126/	213,
145/	132/	216,
160/	126/	237,
224/	150/	285,
232/	150/	309,
256/	214/	349,
271/	204/	376,
277/	216/	388,
280/	216/	397,
307/	216/	412,
349/	240/	460,
352/	240/	469,
337/	276/	472,
352/	276/	517,
391/	324/	544,
400/	336/	565,
496/	370/	685,
512/	402/	717,
544/	386/	733,
520/	402/	741,
592/	482/	829,
656/	450/	861
}
\draw (\a/\c,\b/\c) circle (1/\c)    
;

\draw [->, line width=3.5pt, purple]  (.692+.005,0.538-.012) -- (.692-.005,0.538+.012);    
\draw [->, line width=3.5pt, purple]  (.733-.007, .533-.009) -- (.733+.007, 0.533+.009);   
\draw [->, line width=3.5pt, purple]  (.692+.005, .462+.012) -- (.692-.005, .462-.012);   
\draw [->, line width=3.5pt, purple]  (.75-.01, .5) -- (.75+.01, .5);   
\draw [->, line width=2pt, purple] (.712-.001,.541-.008) -- (0.712+.001,0.541+.006);
\node at (.692-.004,0.538-.017) {$\begin{bmatrix} 4\\ 6 \end{bmatrix}$};
\node at (.733+.0032, .533-.014) {$\begin{bmatrix} 6\\ 3 \end{bmatrix}$};
\node at (.692+0.014, .462+.01) {$\begin{bmatrix} -2\\ \phantom{-}3 \end{bmatrix}$};
\node at (.712-.0005,.541-.017) { \small $\begin{bmatrix} 10\\ 9 \end{bmatrix}$};
\node at (.75-.017,.5-.0055) {  $\begin{bmatrix} 8\\ 0 \end{bmatrix}$};

\end{tikzpicture}
\caption{Apollonian disk packing $(-8,21,24,28)$ and some spinors}
\label{fig:new}
\end{figure}

\noindent
{\bf Example 4:}
Here is  example an  involving a non-symmetric disk packing
generated by a Descartes configuration of disk of curvatures $(-11,21,24,28)$.
The locations of the disk centers may be chosen so that their symbols are integral:
{\small 
$$
\frac{-8,\, -6}{-11}, \quad
\frac{16,\,12}{21}, \quad
\frac{17,\, 12}{24}, \quad
\frac{19,\, 16}{25}
$$}
The resulting spinors are also integral. 
Figure \ref{fig:new} illustrates the packing and a number of spinors
that originate at the disk of curvature $B=24$.
Pick a pair of neighboring spinors  
$$
\mathbf v = \begin{bmatrix} 8\\ 0\end{bmatrix}
\qquad\hbox{and}\qquad
\mathbf w = \begin{bmatrix} \phantom{-}2\\ -3\end{bmatrix}
$$
The corona around that circle can evaluated as
$$
A \ = \ \frac{\pi}{2}  \sum_{\substack {k,n\in \mathbb N \\  \gcd(n,k)=1}}\ 
\frac{1}{\big( 64n^2+32kn+13k^2-24\big)^2}
$$

\section{Comparison with Epstein zeta function}

Recall the definition of the Epstein zeta function \cite{Eps,HW}:
\begin{equation}
\label{eq:epstein}
   Z(S,\rho) = \ \frac{1}{2}\sum_{\mathbf v \in \mathbb Z^n \setminus \{\mathbf 0\}  } 
   \,\left(\,\mathbf v^T\, S \mathbf v\,\right)^{-\rho}, 	
\end{equation}
where $S$ is the $n\times n$ matrix $S$ of a positive definite real quadratic form,  and  $\rho$ is a complex variable
(the real part of which with is greater than $n/2$).
The sum runs over all column vectors with integer entries except the zero vecto. 

The results of the previous section suggest a definition of another arithmetic ``zeta function'', namely:

\begin{definition}
A geometric zeta function is 
\begin{equation}
 Z_\circ (S,\rho,\beta) = \ \frac{1}{2}\sum_{\mathbf v \in \mathbb Z^n_{\Euc}} \,
            \left(\,\mathbf v^T\, S \mathbf v -\beta\,\right)^{-\rho}, 	
\end{equation}
where where $S$ is any $n\times n$ matrix $S$ of a positive definite real quadratic form and  $\rho$ a complex variable.
The sum runs over all column vectors with co-prime integer entries, 
i.e., over the ``$n$=dimensional Euclidean orchard'': 
$$
 \mathbb Z_{\Euc}  \ = \  \{\, \mathbf v = [x_1,..., x_n]^T \in \mathbb Z^n\;|\;    \gcd (x_1,..., x_n) = 1\,\}
$$
Additional condition of geometric compatibility that may be imposed is $\beta = \pm \det S$
\end{definition}

Here is justification for this definition.  First note that, referring to \eqref{eq:main},  
$$
\| M\mathbf u\|^2 =( M\mathbf u)^T\,( M\mathbf u) = \mathbf u^T (M^TM) \mathbf u
$$
where $S = M^TM$ is evidently symmetric positive definite matrix. 
Substituting $S=M^TM$ , $n=2$ and $B=\pm \det M$, and $\rho = 2$
reproduces the formula for the corona area
$$
Area \ = \   \pi \,Z_\circ (M^T\!M,\, 2, \beta) 
\qquad (\, \beta^2 = \pm \det S\,)
$$

Note the difference the two Zeta functions: (1) our sum runs over pairs of co-primes, 
and (2) a shift by $\beta$ is present in the formula.
Analytical properties and possible other arithmetic significance of the new Zeta function remain
to be studied.

\newpage 

\def\coprime{\hbox{\rm coprime}}
\section*{Appendix A:  Practical formulas}

In order to evaluate the area of a corona in a recursive way
one may organize the sum in a diagonal order of counting.
Suppose a disk of curvature $B$ has a pair of neighboring spinors
$$
\mathbf v = \begin{bmatrix} a\\ b\end{bmatrix},
\qquad
\mathbf w = \begin{bmatrix} c\\ d\end{bmatrix},
\qquad 
N = \hbox{some big number}
$$
Then the area between the two disks (end-disk inclusive) is 
$$
\begin{aligned}
K(\mathbf v, \mathbf w) & = \ \sum_{n=1}^N \sum_{k=0}^n \ 
                   \frac{\coprime(n,k)}{ \big( (ka+(n\!-\!k)c)^2 + (kb+(n\!-\!k)d)^2 - B\big)^2}
\\[3pt]
&  = \ \sum_{n=1}^N \sum_{k=0}^n \ \frac{\coprime(n,k)}  
{ \left(\;  \left\| \;\left[\begin{smallmatrix}a&b\\c&d\end{smallmatrix}\right]
                         \left[\begin{smallmatrix}1&-1\\0&\phantom{-}1\end{smallmatrix}\right]
                         \left[\begin{smallmatrix}k\\ m\end{smallmatrix}\right]\;\right\|^2 - B\;\right)^2}
\end{aligned}
$$
where the function $\coprime(n,k)$ returns 1 if $n$ and $k$ are coprimes, and 0 otherwise.

\noindent
Total area of the whole corona may be evaluated by combining the sums over two arcs of the disk:
$$
\begin{array}{rl}
K(\mathbf v,\mathbf w) \ =& \ \dfrac{1}{(a^2 + b^2 -B)^2} + \ \dfrac{1}{(c^2 + d^2 -B)^2} \\[11pt]
                                        & +\sum_{n=2}^N \sum_{k=1}^{n-1} \ 
\left( \dfrac{\coprime(n,k)}{ \big( (am+cn)^2 + (bm+dn)^2 - B\big)^2} 
          + \dfrac{\coprime(n,k)}{ \big( (am-cn)^2 + (bm-dn)^2 - B\big)^2}  \right)
\end{array}
$$
Clearly, the actual value is the limit $N\to\infty$ of the above expressions.

\end{document}